\documentclass{article}
\usepackage[breaklinks=true,colorlinks,citecolor=blue,linkcolor=blue,urlcolor=blue]{hyperref}

\usepackage[margin=3cm]{geometry}
\usepackage{authblk}

\usepackage{natbib}   
\usepackage[latin1]{inputenc}
\usepackage{graphicx}
\usepackage{graphbox}
\usepackage{amssymb}
\usepackage{amsmath}
\usepackage{amsfonts} 	
\usepackage{amsthm} 	
\usepackage{bm} 		
\usepackage{subfig}
\usepackage{multirow}
\usepackage{pbox}
\usepackage{tabularx}
\usepackage{hyperref}
\usepackage{color}
\usepackage{natbib}
\usepackage{verbatim}
\usepackage{listings}
\usepackage{color}
\usepackage{yfonts}
\usepackage{physics}
\usepackage{url}
\usepackage{textcomp}
\usepackage{upquote}

\graphicspath{{./figures/}}

\title{Automatic differentiation for solid mechanics}

\author[1]{Andrea Vigliotti\thanks{andrea.vigliotti@gmail.com, a.vigliotti@cira.it}}
\author[2]{Ferdinando Auricchio}
\affil[1]{ Innovative Materials Laboratory, Italian Aerospace Research Center, 81043 Capua, Italy}
\affil[2]{ Department of Civil Engineering and Architecture, University of Pavia, 27100 Pavia, Italy}

\begin{document}

\maketitle
\date{}

\begin{abstract}
	Automatic differentiation (AD) is an ensemble of techniques that allow to evaluate  accurate numerical derivatives of a mathematical function expressed in a computer programming language.
	In this paper we use AD for stating and solving solid mechanics problems.
	Given a finite element discretization of the domain, we evaluate the free energy of the solid  as the integral of its strain energy density, and we make use of AD for directly obtaining the residual force vector and the tangent stiffness matrix of the problem, as the gradient and the Hessian of the free energy respectively.
	The result is a remarkable simplification in the statement and the solution of complex problems involving non trivial constraints systems and both geometrical and material non linearities.
	Together with the continuum mechanics theoretical basis, and with a description of the specific AD technique adopted, the paper illustrates the solution of a number of solid mechanics problems, with the aim of presenting a convenient numerical implementation approach, made easily available by recent programming languages, to the solid mechanics community.
\end{abstract}

\section{Introduction}

The general problem of solid mechanics hinges on the following familiar pointwise Cauchy's equilibrium equations \citep{Asaro_2006}
\begin{subequations}
	\begin{align}
	\sigma_{ij,j} + b_i &=0 \qquad \text{in} \quad V \label{equil_b}\\
	\sigma_{ij} n_j - t_i &=0 \qquad \text{on}\quad S \label{equil_s},
	\end{align}
\end{subequations}\label{equil}%
which describes the equilibrium of a body subjected to a system of external forces, where $i \in \{1,2,3\}$ are the coordinate directions of a Cartesian frame of reference,  $V$ is the current volume of the body, $S \equiv \partial V$ is the boundary of $V$, $\sigma_{ij}$ are the components of the Cauchy stress tensor, with the letters after the comma in the subscript denoting derivatives along spatial directions and the repeated index denoting summation; furthermore $t_i$ denotes the components of the tractions on $S$, and $b_i$ the components of the body forces in $V$.
Together with stress-strain relationships, and with the kinematic compatibility relations, equations \eqref{equil} define a boundary value problem that allows to find the deformed configuration of a body, given its initial configuration, its material's constitutive laws and a set of boundary conditions.\\
For the simple cases it is possible to find analytical solutions for the above equilibrium equations \citep{Timoshenko1987}.
Nevertheless, in the general practice, solutions are usually sought through numerical methods, such as the Finite Element 
(FE) method  \citep{Zienkiewicz05,Bathe14}.
FE methods are based on a twofold discretization of the problem: (i) the domain and its boundary are discretized in elements connected at nodes; and (ii), the solution is approximated as the weighted sum of a finite set of shape functions, associated with the elements.
The weighting coefficients that control the solution define the Degrees of Freedom (DoFs) of the problem.
In this framework, the problem can be stated and solved, using weighted residual methods such as Galerkin \citep{Zienkiewicz05_SSM}, that are based on selecting the values of the DoFs that minimize a given norm of the error over the domain of interest.
Thus, the problem reduces to finding the solutions of a system of non-linear equations in terms of the DoFs.
The error on each equation defines the residual force vector, and a particular solution can be found  bringing the residual to zero, by means of iterative techniques, such as Newton-Raphson, that make use of the tangent stiffness matrix, which coincides with the Jacobian of the residual vector, to update the trial solution on the basis of its residual.\\

Despite the procedural and algorithmic nature of the general approach, the computer implementation of FE methods for the solution of complex boundary value problems includes many challenging aspects.
With particular reference to solid mechanics (but similar considerations apply to any continuum mechanics problem), the actual statement of equilibrium, \eqref{equil}, involves the modelling of complex material behaviours and requires expressing several vector and tensor quantities, which are best defined in specific, and distinct, reference systems.
The implementation of these steps can in principle be automated, and a number of methodologies have been proposed to facilitate and, at various extents, automate the generation of computer programs capable of efficiently state and solve FE problems.
To this end, the FEniCS project \citep{Logg2007,LoggMardalEtAl2012} set the goal of automating Computational Mathematics Modelling (CMM) problems in general, including the FE methods.
The project aimed at the mechanization of the essential discretization steps of any CMM problem, by means of a suite of general purpose, high level, C++ and python libraries that allow to deal with the numerical implementation of general physical models in a quite abstract, yet efficient, manner, and provide an interface for the definition of the variational problem, its boundary conditions, and its solution.\\

One essential undertaking in the FE implementation of complex solid mechanics problems is translating the mathematical models of the physical processes into a form that can be incorporated into the FE formalism.
This step generally involves analytical manipulations of various sorts, that are normally done by hand and can be source of errors.
In order to address this matter, \cite{Korelc2016} proposed an integrated methodology based on the use of a symbolic engine for the automation of code generation starting from an abstract mathematical statement of the physics under consideration.
In the approach proposed by \citeauthor{Korelc2016} the handling of the symbolic expressions is carried out through the use of \textit{AceGen}, a package within the Mathematica software suite, whose end-product is the source code implementing the models, in different programming languages, for the use in FE programs.
The entire environment includes different components that are capable to produce efficient source code for the generation of the FE residual vector and stiffness matrix.\\

The approach presented here makes use of \textit{Automatic Differentiation} (AD) for the numerical evaluation of the FE residual force vector and stiffness matrix.
In fact, if the material of the body is a \textit{Green elastic} type, for which the deformation work is an exact differential, it is possible to write the expression for the free energy of the solid, and the gradient and the Hessian of the free energy, with respect to the DoFs of the problem, take the meaning of the residual force vector and of the tangent stiffness matrix.
According to the method described in the paper, the gradient and the Hessian of the free energy are not explicitly calculated, but are automatically obtained, through AD, from the function that evaluates the free energy.
The resulting formulation is particularly streamlined and insightful, with the surprising consequence that, with the approach described here, it is possible to write a finite element program without introducing the concept of a stress tensor.
\\

AD is an ensemble of techniques that allows for the numerical evaluation of the derivatives of a function with the same accuracy of the function itself.
AD differs from finite differences, because it does not approximate the continuous derivatives with discrete differences, thus it does not suffer from truncation error, and the only source of error is the inevitable round-off error, due to the finite precision representation of real numbers.
The differentiation techniques based on AD rely on the assumption that the numerical evaluation of the value of a function, regardless of the complexity of the function itself, is always decomposed into in a sequence of elementary sub-expressions by computers.
Therefore, if the analytical derivatives of the sub-expressions are know, it is possible to evaluate the derivatives of the entire function, with respect to the given independent variables, by operating on the partial results.
In this way, AD allows the evaluation of the gradient, along with higher order derivatives, of any computable differentiable function, without having to explicitly write computer code for the derivatives.\\

The paper is structured as follows.
First, we recall the equivalency between the pointwise, strong form statement of equilibrium equations and free energy stationarity.
Subsequently, we discuss the properties of the dual numbers systems,  which is the numerical method chosen to implement automatic differentiation in the present paper.
Lastly we discuss a number of examples of solid mechanics problems that have been solved using AD.
The approach presented in here has been implemented in the Julia programming language, an open source programming language of recent introduction that combines high level of abstraction, high expressiveness and fast execution time \citep{Bezanson2017,Perkel2019}.
All the scripts developed for producing the examples presented in this paper are available under a collaborative licence \citep{GitHubAD4SM}.

\section{Mechanical equilibrium as free energy minimum}\label{sec:free_engy_min}

In conservative systems the deformation work is an exact differential, and it is possible to use the concept of free energy of the system for finding its equilibrium configurations.
In fact, thermodynamics guarantees that all, and only, the configurations that make the free energy stationary are equilibrated \citep{Callen1985}.
Therefore, it is possible to directly write the equilibrium statement by setting the gradient of the system's free energy to nought.
The advantage in using free energy minimization is that free energy is always a scalar quantity, independent of the frame of reference, which is generally well defined and can be calculated using a FE discretization.
In the following we will illustrate the equivalence between free energy minima and equilibria as stated by equation \eqref{equil}, which will allow us to introduce all the relevant solid mechanics quantities as well as to expose the connections between the stress tensor and the gradient of free energy density.\\

Let $\mathfrak{B}$ be a deformable body, occupying a region of an Euclidean space.
Let's assume that $\mathfrak{B}$ is subjected to some external actions from the surrounding environment, in the form of body and surface forces, and in the form of mechanical constraints that restrict its motion.
In its interactions with the environment, the body deforms and can take different configurations.
Let us define the \textit{reference configuration} of $\mathfrak{B}$ as the configuration taken by the body when all of external actions are removed.
Given the reference configuration, the Cartesian coordinates of the points of $\mathfrak{B}$ in any deformed configuration are given by 
\begin{equation}\label{x_i}
x_i = X_i + u_i,
\end{equation}
where $X_i$ are the coordinates of the points of $\mathfrak{B}$ in the reference configuration and $u_i$ are the components of a displacement field mapping the position of the points of $\mathfrak{B}$ from the reference configuration to the \textit{current configuration}.
In association with $u_i$, it is possible to introduce the deformation tensor $\textbf{F}$, as a pointwise measure of the deformation in $\mathfrak{B}$, whose components are given as 
\begin{equation}\label{F_ij}
F_{ij} = x_{i,j} = \delta_{ij} + u_{i,j}
\end{equation}
where $\delta_{ij}$ is the Kronecker symbol and the letter after the comma in the subscript denotes differentiation along direction $i$.
If we assume that $\mathfrak{B}$ is made of a Green elastic material, it is possible to define a strain energy density function, $\phi$, such that the total deformation energy can be expressed as
\begin{equation}\label{Phi_definition}
\Phi =\int_{V_0} \phi\, \mathrm{d}V_0\, ,
\end{equation}
where $V_0$ is the volume occupied by $\mathfrak{B}$ in the reference configuration and $\phi$ is a function of $F_{ij}$ in $V_0$, with units of energy per reference unit volume.
Let us also assume that $\mathfrak{B}$ is subjected to a system of conservative body and surface forces, whose potentials, per unit reference volume and per unit reference surface, are $b_0$ and $t_0$ respectively. 
In accordance with the above assumptions the total free energy of the body is given as
\begin{equation}\label{Psi}
\Psi = \int_{V_0} \left(\phi-b_0\right) \,\mathrm{d}V_0 - \int_{S_0} t_0 \,\mathrm{d}S_0\, ,
\end{equation}
with $S_0\equiv\partial V_0$.
Thermodynamics minimum free energy principle ensures that equilibrium configurations coincide with the configurations that make the free energy stationary.
Therefore, $u_i$ is an equilibrium configuration if and only if
\begin{equation}\label{deltaPsi_0}
\delta \Psi = 0  \qquad \forall\, \delta u_i,
\end{equation}
where $\delta u_i$ is an arbitrary variation in the space of the configurations compatible with the boundary conditions.
We now recall that, through a mere change of variables, the integral \eqref{Psi}, and its variations, can be evaluated in any arbitrary configuration, provided that a mapping exists between the points of the reference configuration and the given configuration. 
Therefore the following holds
\begin{subequations}\label{deltaPsi_1}
	\begin{align}
	\delta \Psi &= \int_{V_0} \left[\frac{\partial \phi}{\partial F_{ij}} \delta F_{ij} - 
	\frac{\partial b_0}{\partial u_i} \delta u_i\right]\, \mathrm{d}V_0 - \int_{S_0}\, \frac{\partial 
		t_0}{\partial u_i} \delta u_i\, \mathrm{d}S_0= \label{deltaPsi_L} \\[5pt]
	&=\int_V \left[J^{-1}\frac{\partial \phi}{\partial F_{ij}} \delta F_{ij} - \frac{\partial 
		b}{\partial u_i} \delta u_i\right]\, \mathrm{d} V - \int_S \frac{\partial t}{\partial u_i} \delta 
	u_i\, \mathrm{d}S \label{deltaPsi_E},
	\end{align}
\end{subequations}
with 
\begin{subequations}\label{deltaPsi_1_sub}
	\begin{align}
	\delta F_{ij} &= \delta u_{i,j}	\label{deltaPsi_1_sub_a}\\[5pt]
	J&=\det\left(\bm{F}\right) \label{deltaPsi_1_sub_b}\\[5pt]
	b&=b_0\,J^{-1} \label{deltaPsi_1_sub_c}\\[5pt]
	t&=t_0\, n_i F_{ik}F_{jk}n_j\, \, J^{-1} \label{deltaPsi_1_sub_d}
	\end{align}
\end{subequations}
where $V$ and $S$ denote the current configuration, as in \eqref{equil}, and $n_i$ are the components of the local normal on $S$.
We observe that with the substitutions \eqref{deltaPsi_1_sub_c} and \eqref{deltaPsi_1_sub_d}, $b$ and $t$ take the meaning of the potential of the external forces per unit current volume and unit current surface, respectively. 
We also observe that through mathematical manipulations, equation \eqref{deltaPsi_E} can be written in terms of displacement variations taken with respect to the current configuration, $\delta \tilde{u}_{i}$, as follows 
\begin{equation}\label{deltaPsi_2}
\begin{split}
\delta \Psi &=\int_V \left[J^{-1}\frac{\partial \phi}{\partial F_{ij}} F_{kj}F^{-1}_{hk} \,\delta 
F_{ih} - \frac{\partial b}{\partial u_i} F^{-1}_{ik}F_{kj}\delta u_j\right]\, \mathrm{d} V -\int_S 
\frac{\partial t}{\partial u_i} F^{-1}_{ik}F_{kj}\delta u_j\, \mathrm{d}S =  \\[5pt]
 &=\int_V \left[\sigma_{ij} \delta\epsilon_{ij} - b_i\, \delta 
\tilde{u}_i\right]\, \mathrm{d} V - \int_S t_i\, \delta 
\tilde{u}_i\, \mathrm{d}S\,, 
\end{split}
\end{equation}
where $F^{-1}_{ij}$ are the components of the inverse of $\bm{F}$, thus $F_{kj}F^{-1}_{hk} = \delta_{jh}$, with the following substitutions:
\begin{subequations}
	\begin{align}
	\sigma_{ij} &= J^{-1}\frac{\partial \phi}{\partial F_{ih}} F_{jh} \label{sigma_ij_a}\\[5pt]
	\delta \tilde{u}_i &= F_{ij}\,\delta u_j \\[5pt]
	\frac{\partial\, \cdot}{\partial \tilde{u}_i} &= \frac{\partial\, \cdot}{\partial u_k} F^{-1}_{ki}  \\[5pt]
	\delta \epsilon_{ij} &= \delta \tilde{u}_{i,j} = F^{-1}_{hj} \,\delta F_{ih} \label{sigma_ij_d}\\[5pt]
	b_i&=\frac{\partial b}{\partial \tilde{u}_i} \\
	t_i&=\frac{\partial t}{\partial \tilde{u}_i}\,,
	\end{align}
\end{subequations} \label{sigma_ij}%
After observing that the following identity holds
\begin{equation}\label{divP}
\left(\sigma_{ij} \delta \tilde{u}_i\right)_{,\,j} = \sigma_{ij,j}\delta \tilde{u}_i + \sigma_{ij} \delta \tilde{u}_{i,j}\,,
\end{equation}
we can express the first variation of the free energy of  $\mathfrak{B}$, making use of the divergence theorem, as
\begin{equation}\label{deltapsi}
\delta \Psi = -\int_{V}  \left(\sigma_{ij,j}+b_i\right) \delta 
\tilde{u}_i \,\mathrm{d}V + \int_{S} \left(\sigma_{ij}n_{j} - t_i\right)\delta \tilde{u}_i  \,\mathrm{d}S = 0 \qquad \forall\, \delta \tilde{u}_i\, .
\end{equation}
Standing the arbitrariness of $\delta \tilde{u}_i$, it follows that each of the integrands in equation \eqref{deltapsi} have to be separately equal to nought everywhere in d$V$ and on d$S$.
Therefore, equation \eqref{deltapsi} is equivalent to equation \eqref{equil}.
We also observe that the same procedure, starting from equation \eqref{deltaPsi_L}, also leads to the equilibrium equation in terms of the nominal stress tensor, or first Piola-Kirchoff tensor, as follows
\begin{subequations}\label{P_ij,j}
	\begin{align}
	P_{ij,j} + b_{0_i} &=0 \qquad \text{in} \quad V_0 \\
	P_{ij}\,N_{j} - t_{0_i}&=0   \qquad \text{on} \quad S_0, 
	\end{align}
\end{subequations}
where $N_j$ are the component of the local normal to the surface on the undeformed configuration, with
\begin{subequations}\label{P_ij}
	\begin{align}
	P_{ij} &= \frac{\partial \phi}{\partial F_{ik}} \label{P_ij_a}\\
	b_{0_i} &= \frac{\partial b_0}{\partial u_i} \label{P_ij_b} \\
	t_{0_i} &= \frac{\partial t_0}{\partial u_i} \,. \label{P_ij_c}
	\end{align}
\end{subequations}
We have thus obtained the equivalence between mechanical equilibrium, in the Newton's laws sense, and the equilibrium in the thermodynamic sense, as the minima of free energy.\\

\section{Automatic differentiation through dual numbers}
In the following we discuss various means for numerically evaluate the derivatives of a function and illustrate the general aspects of automatic differentiation. 
In particular, we will discuss with greater detail the dual number system, which is the frame used to implement AD in the present study.

\subsection{Discrete derivatives approximations}

Finding accurate estimates of the derivatives of multivariate functions at low computational costs is essential in many fields of applied sciences and engineering. 
The simplest numerical estimate for a derivative is \textit{finite difference}.
Finite difference is based on the truncated Taylor expansion formula for analytical functions, and produces the following approximation for a first order derivative
\begin{equation}\label{finite_diff_1st}
\frac{\partial f}{\partial x_i} = \frac{f(\bm{x}+\Delta x_i\,\bm{\imath_{i}})-f(\bm{x})}{\Delta x_i}  + \mathcal{O}\left(\left|\Delta x_i\right|\right),
\end{equation}
where $f$ is a scalar function of the vector $\bm{x}=x_i\bm{i}_i$ while $\bm{i}_i$ are the independent directions of the space of $\bm{x}$ and $\Delta x_i$ is a finite increment for the $i$-th component.
The estimates obtained through equation \eqref{finite_diff_1st} require one additional function evaluation per each independent variable, and suffer from a truncation error of order $\mathcal{O}\left(\left|\Delta x_i\right|\right)$.
We remark that the term $\mathcal{O}\left(\cdot\right)$ does not represents an actual quantification of the approximation of the formula, but it rather represents the convergence rate to the exact value as $\Delta x_i $ collapses.
In addition to truncation error, another important source of inaccuracy is represented by the inevitable round-off, due to the finite precision of the floating point representation of real numbers.
The effects of round-off error are particularly significant in expressions of the type of equation \eqref{finite_diff_1st}, which involves small differences of finite quantities on the numerator, and on the ratio of two small numbers.
As a consequence, it is not possible to arbitrarily enhance the accuracy in the estimate of the derivatives by simply reducing the step length $\Delta x_i$.
One way to improve the accuracy in the  estimation of first derivatives for a fixed $\Delta x_i$ is through central difference scheme as follows
\begin{equation}\label{finite_diff_2nd}
\frac{\partial f}{\partial x_i} = \frac{f(\bm{x}+\Delta x_i\,\bm{\imath_{i}})-f(\bm{x}-\Delta x_i\,\bm{\imath_{i}})}{2\,\Delta x_i}  + \mathcal{O}\left(\left|\Delta x_i\right|^2\right),
\end{equation}
at the cost of two additional function evaluations per variable.
In a similar way, other formulas can be devised that offer improved estimates of the derivative at the cost of a larger number of function evaluation.
We also observe that although the above formulas can be applied recursively for the calculation of higher order derivatives, the accuracy of such estimates rapidly deteriorates since they are based on already noisy estimates of lower order derivatives.\\

\textit{Complex step} is an alternative technique that allows to mitigate round-off error from finite difference \citep{Lyness1967,Lyness1968}.
The complex step method allows to estimate the derivatives of the analytical functions that can be evaluated on the complex plane as follows.
Given the following Taylor series expansion of the function along the imaginary axes of the $i$-th component, 
\begin{equation}\label{complex_step_i}
f(\bm{x}+\imath\,\Delta x_i\,\bm{\imath_{i}}) = f(\bm{x}) +  \frac{\partial f}{\partial x_i} \imath\,\Delta x_i -\frac{\partial^2 f}{\partial x_i\partial x_j} \frac{\Delta x_i\Delta x_j}{2!} + \mathcal{O}\left(\|\Delta \bm{x}\|^2\right),
\end{equation}
where $\imath$ is the imaginary unit, an estimate of the first partial derivative with respect to $x_i$ can be obtained from the imaginary part of the above series as follows
\begin{equation}\label{complex_step_ii}
\frac{\partial f}{\partial x_i}  = \frac{\Im\{f(\bm{x}+\imath\,\Delta x_i\,\bm{\imath_{i}})\}}{\Delta x_i} + \mathcal{O}\left(\left|\Delta x_i\right|^2\right).
\end{equation}
As we can observe, equation \eqref{complex_step_ii} provides an approximation $\mathcal{O}\left(\left|\Delta x_i\right|^2\right)$ of the first derivatives in a single complex valued evaluation of $f(\bm{x})$, which does not suffer from the round-off error due to the difference on the numerator. 
Nonetheless, formulas based on \eqref{complex_step_ii} still need to be evaluated once per each independent variable, at the higher cost of the evaluation on the complex field.\\

\subsection{Symbolic differentiation}

As opposite to finite difference, \textit{symbolic differentiation} treats mathematical expressions as strings of symbols, and applies the analytical differentiation and simplification rules to produce mathematical symbolic expressions that can be evaluated in any programming language.
The availability of robust symbolic differentiation engines has prompted the development of computational approaches tending to the automation of FE code generation \citep{Korelc2016}.
However, the expressions of the derivatives obtained by symbolic differentiation are often far from optimal with respect to computation and memory allocation cost, as they might include redundancy and repetitions that symbolic simplification steps are not capable of removing.
In addition, symbolic differentiation is not directly suitable for dealing with non mathematical functions, like algorithms or computer programs that might include \texttt{for} loops or \texttt{if-then-else} constructs that are common practice in the computer implementation of numerical problems.

\subsection{Automatic Differentiation}

An alternative approach for calculating the numerical values of the derivatives of a function is through \textit{Automatic Differentiation} (AD).
AD is based on the idea that any mathematical expression is evaluated by computers as a sequence of elementary algebraic operations, or call to mathematical functions, with the result being accumulated as the sub-expressions are evaluated \citep{Margossian2018}.
In contrast to finite difference that is based on the discretization of the derivative operator, AD techniques do no try to approximate the differentials but numerically evaluate the derivatives of the sub-expressions, alongside the value of the function itself, using the analytic rules of calculus.\\
Therefore, if the derivatives of all the functions and operators used in the main expression are known, it is possible to evaluate the derivatives of the results with respect to the operands alongside to the value of the function itself.
As a result, AD is not affected by round off errors and allows for the numerical evaluation of  derivatives within the same accuracy of the function itself, with the error only being limited by the machine's representation of floating point numbers.
In addition, AD can be applied recursively to evaluate higher order derivatives with no accuracy detriment, or error build up, because of re-using noisy estimates of lower order derivatives.\\
Two general approaches for the implementation of AD are mostly employed, \textit{forward mode} and \textit{reverse mode}.
In forward mode AD the derivatives of the function, with respect to each of its independent variables, are evaluated along with the main function, in the same order they are encountered, for each sub-expression.
Accordingly, the cost for derivatives evaluation is roughly the cost of one function call per independent variable.
In reverse mode AD the evaluation of the function and of its derivatives takes place in two separate steps.
First, the main function is parsed into sub-expressions that are individually evaluated, whose result is stored, alongside with their derivatives with respect to the argument of each individual sub-expression.
Secondly, the derivatives of the function value with respect to the independent variables are reconstructed using the intermediate derivative of the sub-expressions that have been evaluated, and stored, in the first step.
Thus, at the end of the second step, all the partial derivatives are available in one full run.
Because of its structure, the cost of the reverse mode AD is a few times the cost of the function evaluation alone, depending on how the main expression tree is structured and on how the partial sub-expressions are interconnected, and it is essentially indifferent to the number of the partial derivatives required.
In typical situations we can expect roughly up to a few tens ot times the cost of the evaluation of the function value alone, for a number of partial derivative that can be in the order of more than a few hundred thousands.
Nonetheless, reverse mode AD suffers from the necessity of allocating, and keeping available to the CPU, all of the intermediate results until the entire function is evaluated.
\cite{Griewank2008}, \cite{Hogan2014} and \cite{Elliott2018} provide detailed and thorough description of the forward and the reverse mode AD algorithms.\\

Both forward and reverse mode AD have been traditionally implemented as \textit{algorithmic differentiation} techniques \citep{Bischof_1996,Corliss2002,Naumann2012,Forth2012}, which would take a function's source code as input, and produce the source code for the derivative of the function as output.
While reverse mode AD necessarily requires to operate on the function's source code, forward  mode AD also allows a different type of implementation.
Since forward mode AD is based on a single passage, and it only requires to access the value and the derivatives of the operands at each operation singularly, it can also be implemented by purposely defining a data type that is capable to store both the value and the derivatives of a variable.
Therefore, in the programming languages that allow to extend the ordinary maths operators and function to deal with  user defined data types, through a functionality known as \textit{operator overloading}, once the arithmetic of the extended numerical types is defined, it is possible to call the same code that evaluates a function on the ordinary floating point numeric types, with the newly defined data types, and obtain both the value and the derivatives of the result.
Appendix \ref{sec:implementation_in_Julia} illustrates one application of such technique in the Julia programming language \citep{Bezanson2017,Perkel2019}.\\

In this paper we implement forward mode AD through operator overloading.
Therefore, in the case of conservative systems, it is only necessary to write the code for evaluating the free energy at the integration points of the elements, and the result will also include its gradient and Hessian, that coincide with the element nodal forces and stiffness matrix, respectively.
The choice of operator overloading forces us to use forward mode AD, which is less efficient then reverse mode AD when the number of independent variables increases.
Nonetheless, since we operate on the model element-wise, we only deal with a reduced number of DoFs each time, for instance a QUAD element in 2D involves only 8 DoFs, while a HEXA element in 3D involves 24 DoFs.
At such number of DoFs forward mode is still more efficient than reverse mode because of the reduced no overhead cost needed for preprocessing, and  the absence of additional costs for memory allocation and garbage collection, due to the storage of the intermediate results until the end of the function.\\

\subsubsection{The dual numbers field}\label{sec:dual_number_field}

The particular implementation of the forward mode AD adopted here is based on the recurs to dual numbers, an enriched number system, whose elements have multiple, higher dimensional parts that can deal with the derivative information up to a desired order.
Dual numbers, together with the related algebra, extend real numbers in a way similar to complex numbers and quaternions.
As the complex field allows dealing with expressions that include the square root of negative numbers, and the arithmetic of the quaternions allows to simplify the treatment of rotation in three dimensions, the arithmetic of dual numbers allows the simultaneous calculation of the value of a given expression and its derivatives up to an arbitrary order.

In introducing the dual numbers and their properties, we follow the general treatment of higher 
dimensional number systems as given by \cite{Kantor1989}; a similar treatment of the subject is given by \cite{Fike_2011}.
Howbeit, the treatment of the cited authors did not cover numbers with multiple, separate, higher-order components, with different dimensionality, while the treatment presented in this paper makes use of such structures to deal with derivatives with arbitrary differentiation order.
For the sake of simplicity, in this sections we refer to dual numbers of the second order, and we leave to Appendix \ref{sec:arbitrary_order_dual_numbers} the generalization to dual numbers of arbitrary order. \\

In the present study, dual numbers are the structures chosen to store, and operate, both on the value of a given parameter, $x_0$, and on its derivatives with respect to the independent variables of the problem.
We define dual numbers of the second order and dimension $N$ the quantities of the kind
\begin{equation}\label{mathbfx}
\bm{x} \equiv x_0 + x_i \bm{\imath}_i + x_{ij} \bm{\imath}_{ij} \qquad \text{with} \qquad \begin{aligned}
i &\in 1\dots N \\
j &\in i\dots N
\end{aligned}\,,
\end{equation}
where $x_0$ is the value taken by $\bm{x}$, $x_i$ are the values taken by the first derivatives of $x_0$ with respect to the $i$-th independent variable, and  $x_{ij}$ are the second derivatives of $x_0$ with respect to the $i$-th and the $j$-th independent variable; furthermore, the symbols $\bm{\imath}_i$ are the elements of the canonical base of $\mathcal{R}^N$, while $\bm{\imath}_{ij}$ are defined as
\begin{equation}\label{imath_ij}
\bm{\imath}_{ij} \equiv \bm{\imath}_i\otimes \bm{\imath}_j + \bm{\imath}_j\otimes\bm{\imath}_i\, ,
\end{equation}
with $\bm{\imath}_i\otimes \bm{\imath}_j$ being the dyadic product of $\bm{\imath}_i$ over $\bm{\imath}_j$, thus the following statement holds
\begin{equation}\label{i_ij=i_ji}
\bm{\imath}_{ij} \equiv \bm{\imath}_{ji}\, ,
\end{equation}
which translates the symmetry of the Hessian and,  more in general, the independence of higher derivatives from the order of differentiation.
In the following we refer to the first summand in equation \eqref{mathbfx}, $x_0$, as the real part of $\bm{x}$, and to the second summand, $x_i \bm{\imath}_i$, as its first order dual part and to $x_{ij} 
\bm{\imath}_{ij}$ as the second order dual part or $\bm{x}$.
As we can observe, the first order dual part of $\bm{x}$ is the vector space of real numbers of dimension $N$, while $x_{ij} \bm{\imath}_{ij}$ is the vector space of real symmetric square matrices of dimension N.\\

We will now show that over the set of dual numbers, it is possible to define the operations of addition as well as subtraction, and the operation of multiplication and division.
The neutral element for the sum is the dual zero, i.e. the dual number whose components are all nought, while the neutral element for the product is the dual unity, i.e. the dual number whose value is one, and the derivatives components are all nought.
Therefore, the set of dual numbers as defined above is a field, and the dual zero and unity coincide with the zero and the unity of the real field.\\

By definition, two dual numbers are equal if all of their components are equal, therefore, the following equality statement for dual numbers holds
\begin{equation}\label{mathbfx=mathbfy}
\bm{x} = \bm{y} \iff 
\begin{cases}
y_0 = x_0 \\
y_i = x_i \\
y_{ij} = x_{ij} 
\end{cases}   \qquad \begin{aligned}
\forall\, i &\in 1\dots N \\
\forall\, j &\in i\dots N
\end{aligned}\,.
\end{equation}
For the sake of brevity of notation we will omit to specify the limits value for $i$ and $j$ in the following.\\

We define the sum/the difference of two dual numbers as the sum/the difference of their components, component by component, as follows
\begin{align}\label{z=xpmy}
\bm{z} &= \bm{x} \pm\bm{y} \iff \begin{cases}
z_0 = x_0 \pm y_0 \\
z_i =x_i \pm y_i \\
z_{ij} = x_{ij} \pm y_{ij}
\end{cases}. 
\end{align}
As we can observe, the above definition of sum also induces the definition of the opposite of a dual number as the symmetric with respect to zero for the sum operation.
In a similar way, we define the product of two dual numbers as the sum of all the mixed products of their components, where the following product rules apply for the symbols $\bm{\imath}_i$ and $\bm{\imath}_{ij}$
\begin{subequations}
	\begin{align}
	\bm{\imath}_i \bm{\imath}_j & \equiv \bm{\imath}_{ij}  \label{i_i i_j} \\[5pt]
	\bm{\imath}_{ij} \bm{\imath}_k & \equiv 0 \label{i_ii_ji_k}\, ,
	\end{align}
\end{subequations}
with $\bm{\imath}_{ij}$ defined in equation \eqref{imath_ij}.
With the above positions, the product of two dual numbers follows as
\begin{align}\label{z=xy}
\bm{z} &= \bm{x} \bm{y} \iff \begin{cases}
z_0 = x_0  y_0 \\
z_i = x_iy_0 + x_0y_i\\
z_{ij} = x_{ij}y_0 + x_i y_j + x_j y_i + x_0 y_{ij} 
\end{cases},
\end{align}
we observe that the product of two dual numbers is commutative and associative, thus the following holds
\begin{align}
\bm{xy} &= \bm{yx}  \label{xy=yx} \\
\bm{x}\left(\bm{yz}\right) &= \left(\bm{xy}\right)\bm{z} \label{x(yz)=(xy)z}\, .
\end{align}
As we can observe, given the above definition of multiplication we can define the quotient of two dual numbers $\bm{z}=\bm{x}/\bm{y}$ as the solution to the following equation 
$\bm{zy}=\bm{x}$, therefore the following holds
\begin{equation}\label{z=x/y}
\begin{split}
\bm{z}&=\frac{\bm{x}}{\bm{y}}\, \iff\,\bm{y}\bm{z} = \bm{x} \iff \\ 
	&\iff \begin{cases}
z_0 = \cfrac{x_0}{y_0} \\[10pt]
z_i = \cfrac{x_i}{y_0} - \cfrac{x_0\,y_i}{y_0^2}\\[10pt]
z_{ij} = \cfrac{x_{ij}}{y_0} -\cfrac{x_iy_j+x_jy_i}{y_0^2} + 2x_0\cfrac{y_iy_j}{y_0^3} - 
\cfrac{x_0}{y_0^2}\, y_{ij} 
\end{cases},
\end{split}
\end{equation}
which also defines the inverse of a dual number, as the symmetric to unity with respect to 
multiplication, obtained by replacing $\bm{x}$ with 1 in the above equation.
We observe that, with the definitions given above, the dual numbers field, similarly to the complex number field, is an associative algebra. \\

We now observe that, by recursively applying the identity \eqref{z=xy}, it is possible to extend the operation of raising to integer power over the field of second order dual numbers as
\begin{equation}\label{y=x^n}
\bm{y} = \bm{x}^n \iff \begin{cases}
y_0 = x_0^n \\
y_i = n\,x_0^{n-1}\,x_i \\
y_{ij} = n\left(n-1\right)\,x_0^{n-2}x_ix_j + n\,x_0^{n-1}\,x_{ij} 
\end{cases} .
\end{equation}
More in general, it is possible to extend any continuous, twice differentiable, function $f(x)$ over the second order dual number field by making use of the chain rule for the derivatives as follows
\begin{equation}\label{f(x)}
\bm{f}\left(\bm{x}\right) = f\left(x_0\right) + \frac{\partial f}{\partial x} x_i \,\bm{\imath}_i + \left(\frac{\partial^2 f}{\partial x^2} x_i x_j + \frac{\partial f}{\partial x} x_{ij}\right)\,\bm{\imath}_{ij}\, ,
\end{equation}
where all the derivatives are evaluated in $x=x_0$.
\\

As an example let's assume that we are interested in evaluating the expression 
\begin{equation}\label{y(x1,x2,x3)}
y(x_1, x_2, x_3) = x_1^3\,x_2^2 + x_3^2
\end{equation}
over the second order dual number field, i.e. by treating $x_i$ as second order independent dual quantities, we will show that the result will be a dual quantity itself, retaining the value of the function, and its derivatives with respect to the $x_i$, up to the second order.
Assuming that $x_1$, $x_2$ and $x_3$ are the independent variables, by definition their first derivative with respect to themselves is one, and any other derivatives is zero, therefore their dual representation is the following
\begin{equation}\label{bf{x}_1}
\begin{aligned}
\bm{x}_1 &= x_1 + \bm{\imath}_{1} \\
\bm{x}_2 &= x_2 + \bm{\imath}_{2} \\
\bm{x}_3 &= x_3 + \bm{\imath}_{3}
\end{aligned}\qquad ,
\end{equation}
and the expression \eqref{y(x1,x2,x3)}, evaluated as a dual quantity, takes the following value
\begin{equation}\label{bf{y}i}
\begin{split}
\bm{y}&= (x_1+\bm{\imath}_{1})^3 (x_2+\bm{\imath}_{2})^2 + (x_3+\bm{\imath}_{3})^2 = \\
&= (x_1^3+3x_1^2\bm{\imath}_{1}+3x_1\bm{\imath}_{11}) (x_2^2+2x_2\bm{\imath}_{2}+\bm{\imath}_{22}) + x_3^2 + 2x_3\bm{\imath}_{3} +\bm{\imath}_{33}=\\
&= x_1^3x_2^2+x_3^2 + 3x_1^2x_2^2\,\bm{\imath}_{1} + 
2x_1^3x_2\,\bm{\imath}_{2} + 2x_3\bm{\imath}_{3} +  3x_1x_2^2\,\bm{\imath}_{11} + 6x_1^2x_2\,\bm{\imath}_{12} + x_1^3\,\bm{\imath}_{22} + \bm{\imath}_{33} =\\
&= y_0 + y_i \bm{\imath}_i + y_{ij} \bm{\imath}_{ij}\, ,
\end{split}
\end{equation}
with:
\begin{equation*}\label{bf{y}ii}
\begin{aligned}
y_0 &= x_1^3x_2^2+x_3^2\,, \\
y_i\, \bm{\imath}_i &= 3x_1^2x_2^2\,\bm{\imath}_{1} + 2x_1^3x_2\,\bm{\imath}_{2} + 2x_3\bm{\imath}_{3} \equiv \left[\begin{matrix}
3\,x_1^2x_2^2 \\[5pt]
2\,x_1^3x_2\\[5pt]
2\,x_3
\end{matrix}\right]\,,\\[5pt]
y_{ij}\, \bm{\imath}_{ij} &= 3x_1x_2^2\,\bm{\imath}_{11} + 6x_1^2x_2\,\bm{\imath}_{12} + x_1^3\,\bm{\imath}_{22} + \bm{\imath}_{33} \equiv \left[\begin{matrix}
6\, x_1x_2^2 & 6\,x_1^2x_2 & 0 \\[5pt]
6x_1^2x_2 & 2\,x_1^3 & 0 \\[5pt]
0 & 0 & 2
\end{matrix}\right]\,,
\end{aligned}
\end{equation*}
as we can observe $y_i$ and $y_{ij}$ coincide with the gradient of $y_0$ with respect to $x_i$ and  with its Hessian, respectively.
Appendix \ref{sec:implementation_in_Julia} shows an implementation of the above example in the Julia programming language, carried out with numerical values for $x_1$, $x_2$ and $x_3$.
\\

The dual numbers system can be readily implemented in the programming languages that allow the users to define data types and overload of existing arithmetic operators  over the newly defined types.
The dual numbers type should include data members to hold the real value of the number, and as many higher dimensional arrays up to the desired order of differentiation, with the number of elements in each dimension equal to the number of independent parameters the real part of the number depend on.
Once the data type and the operators have been implemented, scripts with mathematical operations carried out on dual numbers take the same form of the scripts operating on real numbers, and virtually no change is necessary. 
The source code of all the script developed for producing the results presented in this paper are available from the web repositories indicated in \citep{Mendeley}.

\section{Application to solid mechanics}\label{sec:application_to_solid_mechanics}

In this section we will implement AD for the solution of solid mechanics problems.
In first place we will recap how a continuum mechanics problem is generally stated and solved within the FE framework. 
We will then use the same FE discretization for evaluating only the free energy of a solid in a given configuration, and we will recognize how the Jacobian and the Hessian of the system's free energy coincide with the residual force vector and the tangent stiffness matrix, respectively.
Thus, we will show how the complexity of the direct calculation of the residual force vector and the tangent stiffness matrix in FE contrasts with the simplicity of the calculation of the free energy alone.
We will also discuss the implementation, within the AD-assisted framework, of some important components in classical solid mechanics problem, such as non trivial boundary conditions and the hyperelastic material models, while the treatment of geometric non linearity is essentially built-in the AD framework.
\subsection{The Finite Elements framework}
The FE method is generally presented starting from the statement of the principle of virtual work (PVW) \citep{Bathe14,Zienkiewicz05,Bonet08}, which is equivalent to equations \eqref{deltaPsi_1} and \eqref{deltaPsi_2}, in the reference and the current configuration, respectively.
In particular, with reference to \eqref{deltaPsi_L}, the PVW can be expressed, in terms of quantities defined on the reference configuration, as
\begin{equation}\label{VirtualWork}
\int_{V_0} P_{ij}\,\delta F_{ij} \mathrm{d}V_0 =  \int_{V_0} 
\frac{\partial b_{0_k}}{\partial u_k}\, \delta u_k \mathrm{d} V_0 + \int_{S_0} \frac{\partial t_{0_k}}{\partial u_k}\, \delta u_k \mathrm{d}S_0 \quad \forall\, \delta u_k \qquad 
\text{with} \quad \delta F_{ij}=\frac{\partial F_{ij}}{\partial u_k} \,\delta u_k,
\end{equation}
where $\delta u_k$ is a virtual displacement field compliant with the boundary conditions, $\delta F_{ij}$ is the corresponding virtual deformation gradient, and it has been made use of the positions \eqref{P_ij} for the remaining symbols.
The FE approach consists in approximating the evaluation of the integrals in equation \eqref{VirtualWork} by discretizing both the domain of integration and the functional space over which the solution is sought.
The domain of integration is partitioned into elements connected in nodes, and the solution is expressed as the weighted sum of a finite set of shape functions.
The shape functions are defined over the elements, and are selected in a way that guarantees a number of requisites, such as adequate differentiability, continuity across element boundaries and convergence to the analytical solution as the size of the elements collapses.
To this end, the classical choice in FE approaches are multivariate polynomials, which are used to interpolate the nodal values of the unknown function on the interior of each element.
Nonetheless other options are possible, such as the Isogeometric Analysis approach \citep{Hughes2005} where the geometry is described by Non-Uniform Rational B-Splines, and the same rational functions used for the geometry take the role of the shape functions that interpolate the solution over the domain, while the control weights coincide with the coordinates of the vertices of the control polygon for the unknown field.\\

Once the domain has been discretized in elements, and a suitable set  of shape functions has been selected, the integrands in equation \eqref{VirtualWork} are a function of a discrete number of degrees of freedom only.
In this framework the components of the displacement field, and of the displacement gradient at any point of a given element of the domain can be written as
\begin{subequations}\label{N_cdot_u}
	\begin{align}
	u_i &= \bm{N}_i\cdot\bm{u} \\
	F_{ij} &= \bm{N}_{i,j}\cdot\bm{u} + \delta_{ij}
	\end{align}
\end{subequations}
where $\bm{u}$ is the array of the DoFs, $\bm{N}_i$ is the array of the shape functions for the $i$-th component of the displacement field, $\cdot$ is the dot product, and $\bm{N}_{i,j}$ is the $j$-th component of the gradient of $\bm{N}_i$, as usual.
We remark that in this section vector and matrices are marked in bold face, and the product to two vector quantities should be interpreted as the dyadic product, whose result is a matrix.
Therefore, for a given $\bm{u}$, given $\bm{N}_{i,j}$, it is possible to evaluate $u_i$, $F_{ij}$ and any quantity depending on the displacement field, at any point of any element of the FE model of the domain.
In the same way the virtual displacement field can be obtained by means of the same interpolation as
\begin{subequations}\label{N_cdot_delta_u}
	\begin{align}
	\delta u_i &= \bm{N}_i\cdot\delta \bm{u} \\
	\delta F_{ij} &= \frac{\partial F_{ij}}{\partial \bm{u}} \cdot \delta \bm{u}= \bm{N}_{i,j} \cdot \delta \bm{u},
	\end{align}
\end{subequations}
where the arbitrariness of $\delta u_k$ over the functional space of the displacement fields that are compatible with the boundary conditions of the problem translates into the arbitrariness of the components of $\delta \bm{u}$.
Thus, by means of equations \eqref{N_cdot_u} and \eqref{N_cdot_delta_u} it is possible to rewrite equation \eqref{VirtualWork} as
\begin{equation}\label{VirtualWork_DoFs}
\int_{V_0} P_{ij}\frac{\partial F_{ij}}{\partial \bm{u}}\, \mathrm{d}V_0 = \int_{V_0} 
\frac{\partial b_{0}}{\partial \bm{u}}\, \mathrm{d} V_0 + \int_{S_0} \frac{\partial t_{0}}{\partial \bm{u}}\, \mathrm{d}S_0, 
\end{equation}
which is a system of non linear equations, in the unknown unconstrained components of $\bm{u}$.
Therefore, the differential problem of the equilibrium, as stated in equation \eqref{equil}, is translated into a system of non-linear equations, where the unknowns are represented by the unconstrained DoFs.
The integrals in equation \eqref{VirtualWork_DoFs} can be numerically evaluated by means of quadrature rules, and the solution of the FE problem can be found as the zero of the residual vector given by the following
\begin{equation}\label{r_FE}
\bm{r}=\sum_{m=1}^{N_{BE}} \sum_{l=1}^{N_{BW}^m} w_l^m\, \left[P_{ij}\frac{\partial 
	F_{ij}}{\partial \bm{u}} - \frac{\partial b_0}{\partial \bm{u}} \right]_{r^m_l} -
\sum_{m=1}^{N_{SE}} \sum_{l=1}^{N_{SW}^m} v_l^m\, \left[\frac{\partial t_0}{\partial 
	\bm{u}}\right]_{r^m_l}=\bm{0},
\end{equation}
where $N_{BE}$ is the number of volume elements, $N_{BW}^m$ is the number of integration points of the $m$-th volume element, $w_l^m$ is the $l$-th volume integration weight of the $m$-th element; while $N_{SE}$ is the number of surface elements, $v_l^m$ is the $l$-th surface integration weight of the $m$-th surface element, and the subscripts of the square bracket indicate that the quantities enclosed are evaluated at the point ${r^m_l}$, the position of the $w_l^m$ integration weight.

In the finite element formulation equation $\eqref{r_FE}$ can be solved through Newton-Raphson iterative schemes, after an expression for the Jacobian of $\bm{r}$, or the tangent stiffness matrix, has been obtained by differentiating equation $\eqref{r_FE}$ with respect to $\bm{u}$, as follows
\begin{equation}\label{diff(r,u)}
\frac{\partial \bm{r}}{\partial \bm{u}}= \sum_{m=1}^{N_{BE}} \sum_{i=1}^{N_{BW}^m} w_l^m\, 
\left[\frac{\partial P_{ij}}{\partial F_{hk}}\frac{\partial F_{hk}}{\partial \bm{u}}\frac{\partial 
	F_{ij}}{\partial \bm{u}}  - 
\frac{\partial^2 b_0}{\partial \bm{u} \partial \bm{u}}\right]_{r^m_l} - \sum_{m=1}^{N_{SE}} 
\sum_{i=1}^{N_{SW}^m} v_l^m\, \left[\frac{\partial^2 t_0}{\partial \bm{u}\partial 
	\bm{u}}\right]_{r^m_l} \quad , 
\end{equation}
where we made us of the fact that, since $F_{ij}$ is linear in $\bm{u}$, the following holds
\begin{equation}\label{d2Fdudu}
\frac{\partial^2 F_{ij}}{\partial \bm{u} \partial \bm{u}} = \bm{0}\, .
\end{equation}
The calculation of the summands in equation \eqref{r_FE} and \eqref{diff(r,u)} is the core of the FE methods and represent the most challenging aspect of the computer implementation of the method.
In particular the calculation the components of the stress tensor, $P_{ij}$, and their derivatives ${\partial P_{ij}}/{\partial F_{hk}}$, is in general a quite sophisticated task, since it requires dealing with second order and fourth order tensors respectively.

\subsection{The automatic differentiation formulation}

We now turn our attention to the use of the AD for the solution of the equilibrium problem.
We begin by assuming that the material of the solid is a Green elastic material, because under this assumption the resulting formulation is particularly simple and insightful.
For Green elastic solid the deformation work is an exact differential, and the free energy function for the system is given by equation \eqref{Psi}.
In order to numerically evaluate the integrals in equation \eqref{Psi} we can make use of the same twofold discretization used for FE, obtaining the following expression
\begin{equation}\label{Psi_FE}
\Psi(\bm{u}) = \sum_{m=1}^{N_{BE}} \sum_{i=1}^{N_{BW}^m} w_l^m \left[\phi + b_{0}\right]_{r^m_l} + 
\sum_{m=1}^{N_{SE}} \sum_{i=1}^{N_{SW}^m} v_l^m\, \left[t_0\right]_{r^m_l}\,,
\end{equation}
where it has been highlighted that the free energy, within the FE discretization, is a function of the array of DoFs, $\bm{u}$.
As discussed in section \ref{sec:free_engy_min}, equilibrium configurations are those that satisfy equation \eqref{deltaPsi_0}, which in the FE discretization can be written as
\begin{equation}\label{deltaPsi_FE}
\delta \Psi= \frac{\partial\, \Psi}{\partial \bm{u}} \cdot \delta\bm{u} = 0\,, \, \forall \,\delta\bm{u} \, \iff \, \frac{\partial\, \Psi}{\partial \bm{u}}=\bm{0}\,.  \hspace{20mm}
\end{equation}
The expression above is a system of non-linear equations in the unknown $\bm{u}$, whose residual and Jacobian, are given as the gradient and Hessian of $\Psi$, respectively, as
\begin{subequations}
	\begin{align}
	\bm{r} &= \frac{\partial\, \Psi}{\partial \bm{u}} \label{dPhi_du} \\[5pt] 
	\frac{\partial \bm{r}}{\partial \bm{u}} &= \frac{\partial^2\, \Psi}{\partial \bm{u} \partial 
		\bm{u}}\,. \label{dPhi_dudu}
	\end{align}
\end{subequations}
We remark that the residual on equation \eqref{r_FE}, which derives from the PVW statement given in \eqref{VirtualWork}, and the residual on equation \eqref{deltaPsi_FE}, which is obtained as the first the deformation work in the equilibrium configuration given in \eqref{Psi_FE}, are both work-conjugated through the same virtual nodal displacements, $\delta\bm{u}$, hence they must coincide.
Therefore the gradient of the free energy coincides with the residual force vector of the finite element problem.
At the same time, the tangent stiffness matrix coincides with the Hessian of the free energy, being both the derivative of $\bm{r}$ with respect to $\bm{u}$.\\

We now remark that the gradient and the Hessian of $\Psi$, in equations \eqref{dPhi_du} and \eqref{dPhi_dudu}, respectively, can be both immediately calculated by the same computer program that evaluates equation \eqref{Psi_FE}, through automatic differentiation, if $\bm{u}$ is treated as an array of dual numbers, and the dual number algebra has been implemented in the programming language.
Therefore, AD allows to evaluate the residual force vector and the tangent stiffness matrix by simply calculating the numerical integral of the free energy density over the domain.\\

We finally remark that the same approach can be used in the cases when a functional relation exists between $P_{ij}$ and $F_{ij}$, but no elastic potential can be defined.
For these materials equation \eqref{dPhi_dudu} is replaced by \eqref{r_FE}, through equations \eqref{N_cdot_u} and \eqref{N_cdot_delta_u}, with the components of $\bm{u}$ being independent dual quantities, while equation \eqref{dPhi_dudu} still holds and it is obtained as the first order dual components of $\bm{r}$. \\

In the sections that follow we will illustrate how some of the fundamental elements in a solid mechanics problem, such as non linear constitutive laws, or complex boundary conditions, can be easily included in the problem formulation with recurs to automatic differentiation for their implementation.

\subsection{Boundary conditions and constraint equations}\label{sec:constraints}

Non trivial boundary conditions can be applied with the use of Lagrange multipliers using automatic differentiation technique for the direct evaluation of the gradient and of the Hessian of the Lagrange function.
We recall that, following the Lagrange multipliers technique, the minimization of a function, in the presence of constraints can be achieved by weighting the residuals of the constraint equations through unknown factors, the Lagrange multipliers, and adding them to the function to be minimized, as follows
\begin{equation}\label{L}
L\left(\bm{u}, \bm{\lambda}\right) = \Psi\left(\bm{u}\right) - \bm{\lambda}\cdot 
\bm{g}\left(\bm{u}\right),
\end{equation}
where $\Psi$ is the function to be minimized in the first place, which in our case is the free energy of the solid, $\bm{u}$ are the degrees of freedom of the problem, $\bm{g}\left(x_i\right)$ is the array of the constraint equations and $\bm{\lambda}$ is the array of the Lagrange Multipliers.

Since the problem of minimizing $L$ with respect to $\bm{u}$ and $\bm{\lambda}$, is essentially identical to the problem of minimizing $\Psi$ with respect to $\bm{u}$ only, it can be treated in the same way.
Nonetheless, since $L\left(\bm{u}, \bm{\lambda}\right)$ is linear in $\bm{\lambda}$, it is not necessary to treat $\lambda$ as dual quantity, but suffices to evaluate $\Psi\left(\bm{u}_i\right)$ and $\bm{g}\left(\bm{u}_i\right)$ over the dual number field of $\bm{u}$, while the components of the augmented gradient and Hessian of $L$ can be obtained from the expression of the first and second variation of $L$, which is given as
\begin{subequations}\label{deltaL_delta^2L}
	\begin{align}
	\delta L &= \left(\frac{\partial \Psi}{\partial \bm{u}} - \bm{\lambda} \cdot \frac{\partial \bm{g}}{\partial \bm{u}} \right)\cdot \delta \bm{u} -\bm{g} \cdot \delta \bm{\lambda} \\
	\begin{split}
	\delta^2 L &= \left(\frac{\partial^2 \Psi}{\partial \bm{u}\partial \bm{u}} - \bm{\lambda} \cdot \frac{\partial^2 \bm{g}}{\partial \bm{u}\partial \bm{u}}\right)\colon \delta \bm{u} \delta \bm{u} - \frac{\partial \bm{g}}{\partial \bm{u}} \colon \delta \bm{u} \delta \bm{\lambda} + \\
	 &\qquad - \delta \bm{\lambda}\cdot\,\frac{\partial \bm{g}}{\partial\bm{u}}\cdot\delta \bm{u} \,,
	\end{split}
	\end{align}
\end{subequations}
which yields the following, in block matrix notation,  
\begin{subequations}\label{deltaL_delta^2L_matrix}
	\begin{align}
	\nabla L &= \left[\begin{matrix}
	\dfrac{\partial \Psi}{\partial \bm{u}} - \bm{\lambda} \cdot \dfrac{\partial \bm{g}}{\partial \bm{u}}\\[5pt]
	-\bm{g}
	\end{matrix}\right] \label{deltaL}\\[5pt]
	\nabla^2 L &= \left[\begin{matrix}
	\dfrac{\partial^2 \Psi}{\partial \bm{u}\partial \bm{u}} - \bm{\lambda} \cdot \dfrac{\partial^2 \bm{g}}{\partial \bm{u}\partial \bm{u}} & \hspace{20pt}-\dfrac{\partial \bm{g}}{\partial \bm{u}}^T \\[5pt]
	-\dfrac{\partial \bm{g}}{\partial \bm{u}} & \hspace{20pt}\bm{0}
	\end{matrix}\right]	\,. \label{delta^2L}
	\end{align} 
\end{subequations}
Therefore the problem reduces to solving the following
\begin{equation}\label{nablaL=0}
\nabla L = 0 \,,
\end{equation}
where equation \eqref{delta^2L} takes the meaning of the tangent stiffness matrix of the problem.

\subsection{Hyperelastic material models}\label{sec:materials}
We now turn our attention to the most common expressions for the strain energy density functions of materials.
In very simple cases, such as for the small deformations of strut or beam elements, the stress is uniaxial, the deformation state of the solid is adequately described by a single scalar quantity, and the material behaviour can be treated as linear elastic, with a strain energy density function of the type
\begin{equation}\label{Phi_Hooke}
\Phi^H = \frac{1}{2} E_s \left(1+\epsilon_n\right)^2,
\end{equation}
where $E_s$ is the Young modulus of the material, and $\epsilon_n$ the component of the nominal strain tensor conjugated to only non zero stress tensor component .
In these cases the calculation of the deformation work is particularly simple; however, for a general solid mechanics problem, the state of deformation has arbitrary principal directions and distinct principal stretches, therefore, more sophisticated expressions for the deformation energy density are used.
Green elastic materials are a quite general class of material models for which the strain energy density function is assumed as a local function of the components of the right Cauchy-Green deformation tensor, $\bm{C}=\bm{F}^T\bm{F}$ \citep{Ogden_2013}, as
\begin{equation}\label{Phi_Green}
\phi^\text{G}=\phi^\text{G}\left(C_{ij}\right),
\end{equation}
where $C_{ij}$ are the components of  $\bm{C}$.
Among Green elastic materials, one class of material models that are of great interest for engineering applications are the isotropic hyperelastic materials, whose behaviour is invariant to rigid rotations of the applied strain.
For these material models, the strain energy density function can be expressed as a function the invariants of $\bm{C}$ only, defined as
\begin{equation}\label{I_i}
\begin{split}
& \qquad \qquad \phi^\text{iso}= \phi^\text{iso}\left(I_1, I_2, I_3\right) \quad \text{with:} \\
& \begin{aligned}
I_1 &= C_{11}+C_{22}+C_{33}\\
I_2 &= C_{11}C_{22}+C_{22}C_{33}+C_{11}C_{33}-C_{21}^2-C_{31}^2-C_{32}^2\\
I_3 &= C_{11}C_{22}C_{33}+2C_{21}C_{31}C_{32}-C_{11}C_{32}^2-C_{22}C_{31}^2-C_{33}C_{21}^2
\end{aligned}
\end{split}
\end{equation}

Incompressible isotropic materials are subject to the internal isochoric constraint, $I_3=1$, therefore $\phi^\text{iso}$ depend on $I_1$ and $I_2$ only, and its general expression is of the type \citep{Ogden_2013}
\begin{equation}\label{phi_inc}
\phi^\text{inc} = \sum_{p,q=0}^{\infty} c_{pq}\left(I_1-3\right)^p\left(I_2-3\right)^q \qquad 
\text{with} \quad I_3=1 \,,
\end{equation}
where $p$ and $q$ are non negative integers and $c_{pq}$ are non negative real parameters.
Expression \eqref{phi_inc} is only valid if the isochoric constraint is explicitly enforced, however, since such a constraint often yields to convergence problems  in the finite element formulation, when using an incompressible material model, the following decomposition of the $\bm{C}$ can be assumed,
\begin{equation}\label{C=CvCiso}
\bm{C}=\overline{\bm{C}}\bm{C}^\text{v} \qquad \text{where}\qquad \begin{aligned}
\bm{C}^\text{v} &= J^{2/3}\textbf{I}\\
\overline{\bm{C}} &= J^{-2/3}\bm{C}
\end{aligned}\,,
\end{equation}
where $J=\det\left(\bm{F}\right)$ and $\textbf{I}$ is the identity tensor.
It is straightforward to verify that $\det\left(\overline{\bm{C}}\right)=1$, therefore equation \eqref{C=CvCiso} decompose the total deformation into a isochoric deformation, represented by $\overline{\bm{C}}$, and a purely volumetric deformation, given by $\bm{C}^\text{v}$.
Following the decomposition \eqref{C=CvCiso}, $\Phi^\text{inc}$ is approximated as
\begin{equation}\label{Phi}
\Phi = \sum_{p,q=0}^{\infty} c_{pq}\left(\overline{I}_1-3\right)^p\left(\overline{I}_2-3\right)^q + f\left(J\right)\,,
\end{equation}
where $\overline{I}$ are the invariants of $\overline{\bm{C}}$, and $f\left(J\right)$ is a positive function of $J$ that effectively penalizes volume variations.
In the example section of this paper we will use, in particular, the Mooney-Rivlin and the Neo-Hokkean models, whose strain energy density expression is given by 
\begin{subequations}\label{hyperelas}
	\begin{align}
\begin{split}
	\text{\small Neo-Hookean} 	\\
	 \quad \phi^{\text{NH}}&=c_{10}(\overline{I}_1-3) + G \left(J-1\right)^2   
\end{split} \label{NeoHookean} \\ 
\begin{split}
	\text{\small Mooney-Rivlin} \\
	\quad \phi^{\text{MR}}&=c_{10}(\overline{I}_1-3) + c_{01}(\overline{I}_2-3) + G \left(J-1\right)^2
\end{split} \label{MooneyRivlin}
	\end{align}
\end{subequations}
where $c_{10}$, $c_{01}$ and $\text{G}$ are constant, non negative, parameters that define the material behaviour.\\

\subsection{Derivation of the stress tensor in an AD framework}
We remark that in the approach presented here, since the residual force vector and the tangent stiffness matrix are automatically obtained from the free energy, we never explicitly calculate the components of the stress tensor.
However, the value of the entries of the stress tensor are still important quantities, since resistance criterion, such as Von Mises, are based on it.
Nonetheless, they can always be evaluated, as a post processing step, from the equalities \eqref{sigma_ij_a} and \eqref{P_ij_a}, by means of the applicable expression for the strain energy density, by treating the components of $\bm{F}$ as independent dual quantities, whose value is obtained from the displacement field of the equilibrium configuration.

\section{Examples}\label{sec:Examples}

In this section we present a selection of solid mechanics problems whose solution has been found with the recurs to the automatic differentiation techniques described in the paper.
The examples presented include structural elements, such as rods and beams (section \ref{sec:truss} and \ref{sec:Euler_beams}), continuous plane stress elements (section \ref{sec:plane_stress}), a problem with cylindrical symmetry (section \ref{sec:cylindric_symmetry}), and a full three-dimensional problem (section \ref{sec:3D_solid}).\\
All the problems included the effects of geometric non-linearities, non trivial boundary conditions, and the hyperelastic material models described in section \ref{sec:materials}.
The non trivial boundary conditions were introduced using the Lagrange multipliers technique, as described in section \ref{sec:constraints}.\\
All problems presented here were solved using the Julia programming language \citep{Bezanson2017,Perkel2019}, and the script files used for the solution have been made available to the reader \citep{Mendeley}.

\subsection{The non linear truss}\label{sec:truss}

\begin{figure*}
	\centering
	\includegraphics[width=0.9\textwidth]{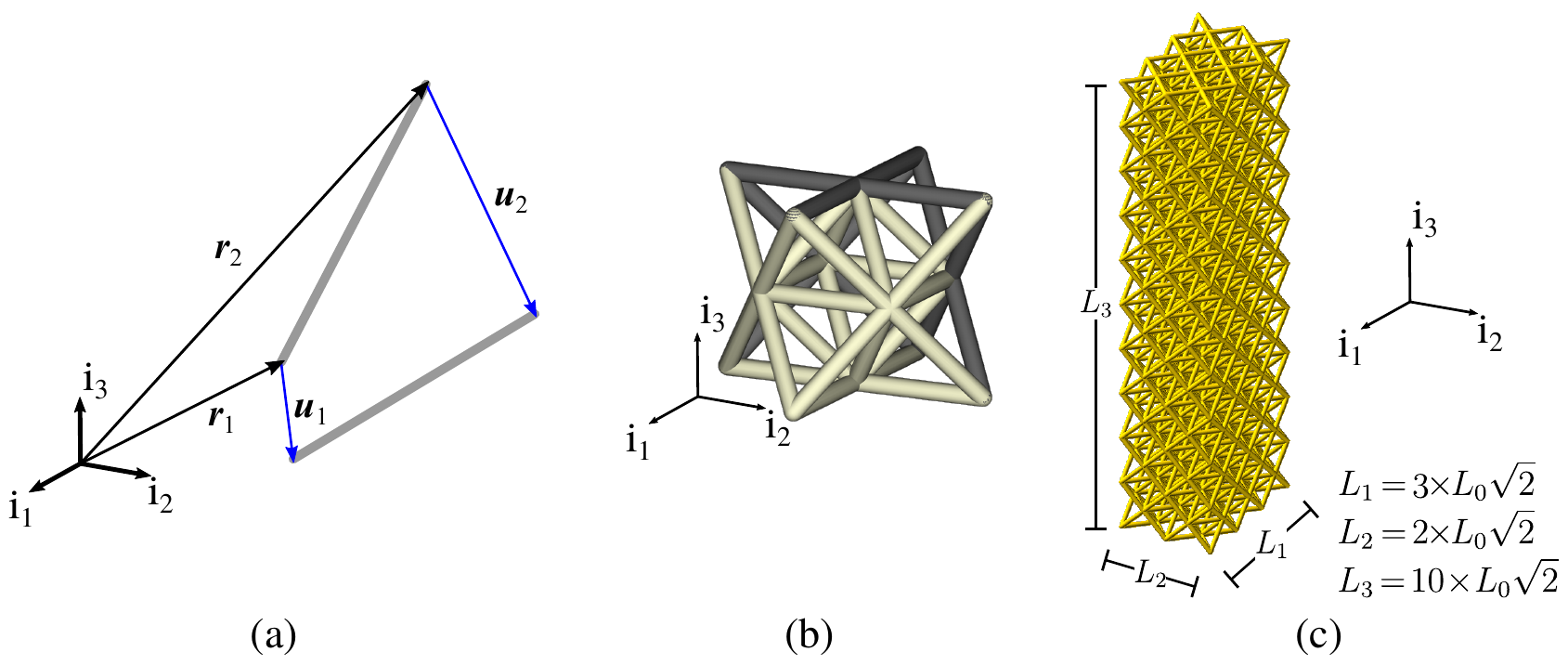}
	\caption{(a) Rod element, $\bm{r}_1$ and $\bm{r}_2$ are the nodes positions in the reference 
		configuration, $\bm{u}_1$ and $\bm{u}_2$ are the displacement vectors. (b) Sketch of the 
		regular octet unit cell, all struts have the same length $L_0$.(c) Sketch of the entire truss structure}\label{fig:fig_rods_UC_truss}
\end{figure*}

In this section we consider the equilibrium of a tridimensional structure made of prismatic elements, connected at their endpoints to form a truss.
We also assume that the material of the struts is linear elastic, with Young modulus $E_s$, and  that cross section deformations are negligible with respect to the axial deformation of the elements.
Under these  assumptions, the elements can only store elastic energy by variations of their length, and the deformation energy of a single element, is given as
\begin{equation}\label{phi_rods}
\phi^{\text{rod}} = A\, l_0 \Phi^H,
\end{equation}
where $A$ is the cross section area, $l_0$ is the reference length, $\Phi^H$ is defined by equation \eqref{Phi_Hooke}, with $\epsilon_n = l/l_0-1$, and $l$ is the length of the element in the current configuration.
The total deformation energy of the truss can then be readily obtained as the sum of the deformation energy of all of its elements, and it is given as
\begin{equation}\label{Phi_truss}
\Phi^{\text{truss}} = \sum_i \phi_i^{\text{rod}},
\end{equation}
where $\phi_i^{\text{rod}}$ is the strain energy of the $i-$th element.
With reference to figure \ref{fig:fig_rods_UC_truss}.a, let $\bm{r}_1$ and $\bm{r}_2$ be the positions of the end nodes of a rod in the reference configuration, and $\bm{u}_1$ and $\bm{u}_2$ be the displacement vectors of the nodes, $l$ and $l_0$ can be easily obtained for a given element as
\begin{equation}\label{l0_and_l}
\begin{aligned}
l_0 &= \|\bm{r}_2-\bm{r}_1\| \\
l   &= \|\bm{r}_2+\bm{u}_2-\left(\bm{r}_1+\bm{u}_1\right)\| .
\end{aligned}
\end{equation}
Therefore, since a truss can be idealized as network of rods connecting in nodes, given the topology of the connections, and the cross section and material properties of the struts, the deformation energy of a truss can be easily computed as function of the components of the displacements of the nodes.\\

Here we consider the equilibrium of a structure obtained  by replicating the regular octet unit cell, shown in figure \ref{fig:fig_rods_UC_truss}.b along the directions $\bm{\imath}_{1}$, $\bm{\imath}_{2}$ and $\bm{\imath}_{3}$, without duplicating the coincident rods.
The regular octet is a well known structure, which is characterized for its lightness and strength \citep{Fuller1966,Deshpande01b}.
Since the regular octet topology is both statically and kinematically determined, no mechanisms arise from its deformation, and it can withstand any external load by stretching of its elements only \citep{Deshpande01}.
In particular we consider a structure made of 2 unit cells in the directions $\bm{\imath}_{1}$, $N_1=2$, three units in direction $\bm{\imath}_{2}$, $N_2=3$, and 10 cells in direction $\bm{\imath}_{3}$, $N_3=10$, as shown in figure \ref{fig:fig_rods_UC_truss}.c.
\\

\begin{figure*}
	\centering
	\includegraphics[width=0.9\textwidth]{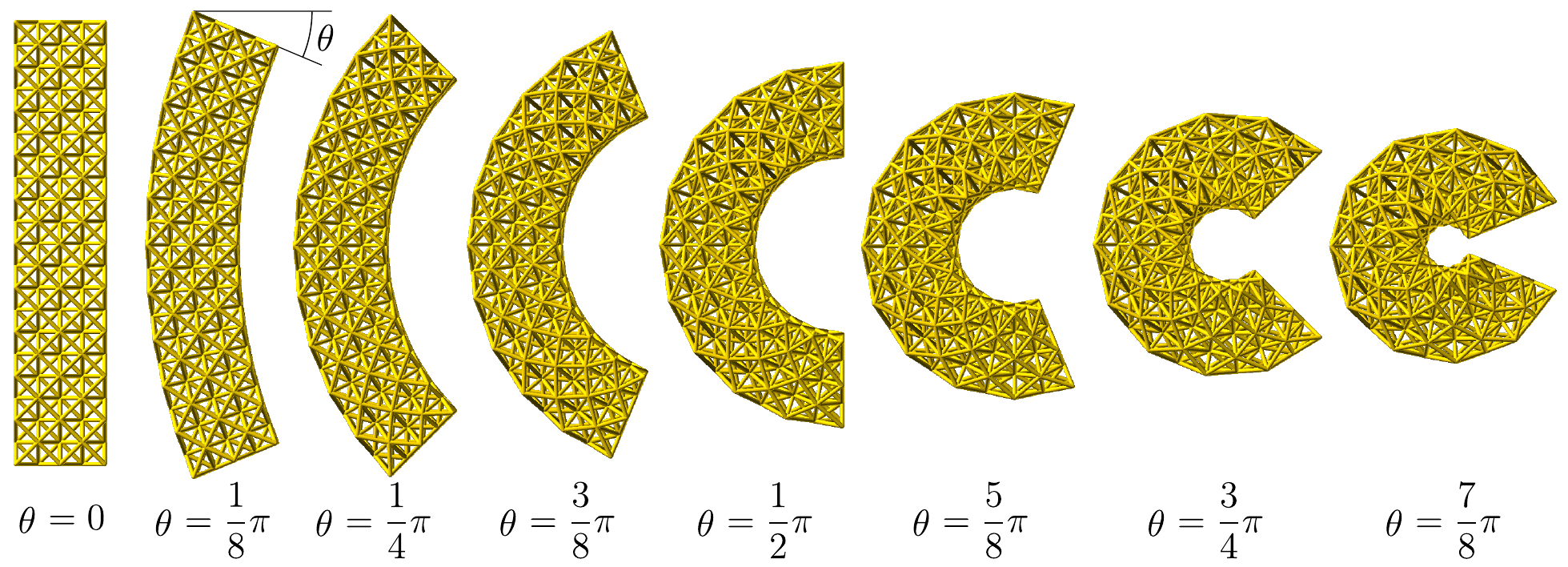}
	\caption{Shape of the lattice structure for different rotation angle of the top and bottom planes. The rotation on the top and the bottom planes have equal amplitude and opposite sign.}\label{fig:nodes_01_def_casethree}
\end{figure*}

The selected boundary conditions produce the bending of the domain around an axis parallel to $\bm{\imath}_{1}$ by constraining the nodes on both the the ends of the domain to remain on planes that are symmetrically rotated around an axes parallel to $\bm{\imath}_{1}$, for a prescribed angle $\theta$.
Under such boundary conditions the rigid body translation along direction $\bm{\imath}_{1}$ is still available to the structure, and it is removed by constraining any motion of the center of gravity in direction $\bm{\imath}_{1}$.
In summary, the boundary conditions are implemented through Lagrange multipliers and are expressed as
\begin{subequations}
	\begin{align}
	\left(X_2+u_2\right) \cos\left(\theta\right) + \left(X_3+u_3\right) \sin\left(\theta\right) - \frac{1}{2}L_3 &= 0 \label{C_eq_truss_a}\\[5pt]
	\left(X_2+u_2\right) \cos\left(\theta\right) - \left(X_3+u_3\right) \sin\left(\theta\right) + \frac{1}{2}L_3 &= 0 \label{C_eq_truss_b}\\[5pt]
	\sum u_1 = 0 \,,\hspace{80pt} \label{C_eq_truss_c}
	\end{align}
\end{subequations}
where the condition \eqref{C_eq_truss_a} applies to the nodes at $X_3=L_3/2$, condition \eqref{C_eq_truss_b} applies to the nodes at $Z_3=-L_3/2$, while the summation in \eqref{C_eq_truss_c} extends to all the nodes of the model and serves the purpose of removing the residual rigid body degrees of freedom.
The model comprised 368 nodes and 2160 rods, boundary conditions were applied in $\pi/8$ steps between $\pi/8$ and $7\pi/8$.
In all cases convergence was achieved in 5 iterations, except in steps 6 and 7 where convergence was achieved in 6 and 7 steps respectively, which confirmed that the Hessian of the free energy was accurate.
Simulation results for different values of $\theta$ between 0 and $7\pi/8$ are shown in figure \ref{fig:nodes_01_def_casethree}.
As it can be observed, the formulation adopted was capable to attain convergence also in the presence of very large displacements, with each step taking only three iterations to converge.

\subsection{Euler beams}\label{sec:Euler_beams}

\begin{figure}[h!]
	\centering
	\includegraphics[width=60mm]{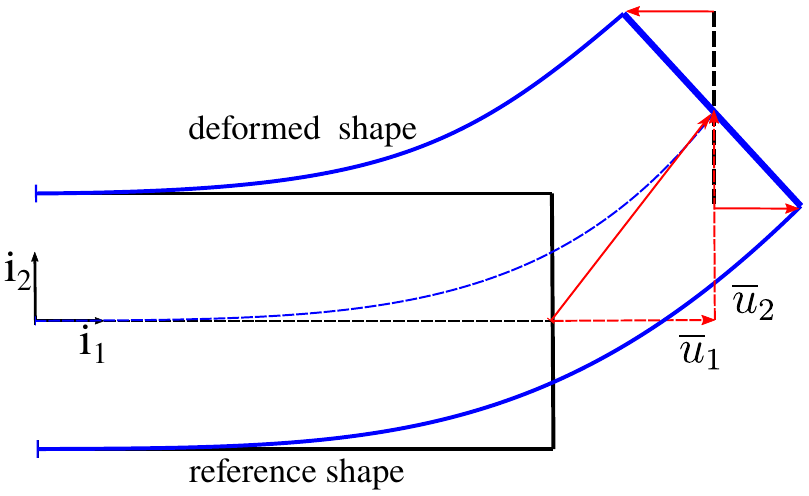} 
	\caption{Euler beam displacement model}\label{fig:fig_beam01}
\end{figure}

In a similar way it is possible to analyse the response of structures made of Euler beams.
Given a prismatic bar, the Euler beam model assumes that each cross section rigidly rotates around an axis orthogonal to the beam axis, passing through the centre of gravity of the section, neglecting any shear contribution to deformation energy and load bearing.
With reference to figure \ref{fig:fig_beam01}, in a reference frame with $\bm{\imath}_{1}$ aligned with the beam axis, and $\bm{\imath}_{2}$ aligned with the axis of rotation of the cross section, under the assumption of small local cross section rotation, the following displacement model for the points of the beam holds
\begin{equation}\label{uv_beams}
\begin{aligned}
u_1 &= \overline{u}_1 - X_2\,\overline{u}_{2,1}  \\[5pt]
u_2 &= \overline{u}_2,
\end{aligned}
\end{equation}
where $\overline{u}_i$ are the components of the displacement of the points on the beam axis,  $X_i$ are the coordinate of the points of the beam in the reference configuration, and the subscript after the comma denotes differentiation with respect to the coordinate $X_i$.
Since the shear contributions to deformation energy are neglected, under finite displacement assumption, the elastic energy is only a function of the first component of the Green Lagrange deformation tensor, which, in accordance to equation \eqref{uv_beams}, is given as
\begin{equation}\label{E11_beam}
\begin{aligned}
E_{11} &= \frac{1}{2} \left[\left(1+u_{1,1}\right)^2+u_{2,1}^2-1\right] = \\
&= \frac{1}{2} \left[\left(1+\overline{u}_{1,1}-X_2\,\overline{u}_{2,11}\right)^2 + \overline{u}_{2,1}^2-1\right].
\end{aligned}\end{equation}
Under the assumption that a linear elastic model is adequate for representing the material behaviour, the deformation energy of an individual beam and of a structure made of beam elements are given respectively as
\begin{align}
\phi_i &= \frac{1}{2} E_s \int_{V_i} E_{11}^2 \,\mathrm{d} V \label{phi_beam} \\[5pt]
\Phi &= \sum_i \phi_i \label{Phi_beam}  .
\end{align}
Therefore, once a parametric representation for the displacement of the points of the axis of the beams is given, the deformation energy of the structure can be expressed as a function of the chosen parameter, and the equilibrium configuration can be found by minimizing the free energy of the structure.
A common parametric representation for $\bm{\overline{u}}$ assumes the nodal displacement and rotations as the degrees of freedom of the elements, and takes the axial component of the displacement as a linear function of $X_1$, and the transverse components as cubic functions of $X_1$ \citep{Bathe14}.\\

\begin{figure}[h!]
	\centering
	\subfloat[]{\includegraphics[height=5cm]{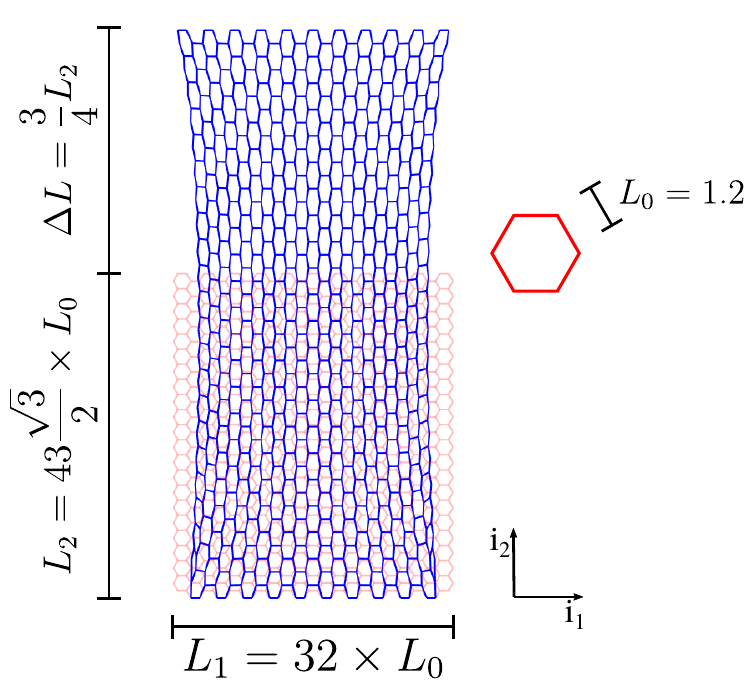} \label{fig:HexaLattice_undef_def}} \quad
	\subfloat[]{\includegraphics[height=5cm]{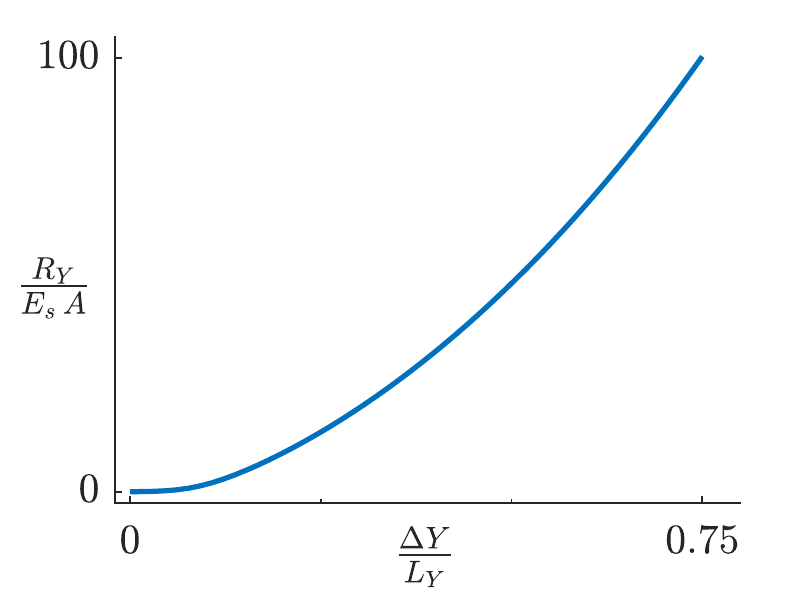}\label{fig:HexaLattice_02}}
	\caption{Simulation results for a bidimensional hexagonal lattice made of Euler beam elements, 
		the beam elements have square cross section with side $t=L_0/10$. (a) Deformed and undeformed 
		configuration of the lattice for total applied. (b) Normalized total reaction force vs applied 
		displacement curve, $R_2$ is the sum of the residual conjugated to the vertical component of 
		the displacement of the top nodes, $E_s$ is the Young Modulus of the material, 
		$A=\left(L_0/10\right)^2$ is the cross section area of the elements.}\label{fig:HexaLattice}
\end{figure}

Figure \ref{fig:HexaLattice} shows the equilibrium configuration of a bidimensional hexagonal lattice under tension obtained using equation \eqref{Phi_beam} for evaluating the deformation energy and the automatic differentiation approach described in this paper to find the stationary energy configurations.
Displacement boundary conditions have been applied to the nodes on the top and bottom of the model, by preventing horizontal displacement and rotation around the axis orthogonal the plane of the model, and prescribing the vertical displacement of the top nodes.
All elements have the same length, $L_0$, and the same square cross section with side $t= L_0/10$.
Figure \ref{fig:HexaLattice_undef_def} shows the deformed configuration for a prescribed total displacement in the vertical direction of $\Delta L=0.75\, L_2$, where $L_2$ is the initial length of the model in the direction $\bm{\imath}_{2}$, while figure \ref{fig:HexaLattice_02} shows the total reaction force, obtained as the sum of the residuals conjugated to vertical displacement for the top nodes of the model,  normalized by the Young Modulus of the material, $E_s$, and the cross section area, $A=t^2$, as a function of the applied displacement.
The solution was obtained in 20 steps, with every step taking between 6 and 9 Newton-Raphson iterations.

\subsection{The plane stress problem}\label{sec:plane_stress}

\begin{figure}[h!]
	\centering
	\includegraphics[width=0.5\textwidth]{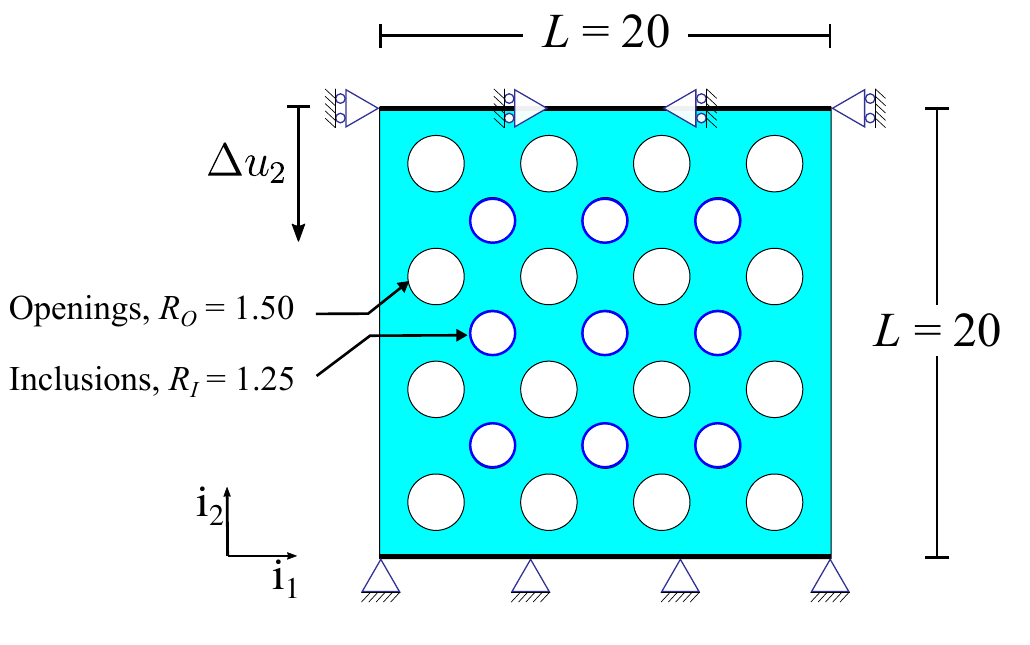}
	\caption{Plane stress problem, undeformed domain with boundary conditions, $\Delta u_2$ is the 
		applied displacement, the radius of the openings is $R_O=1.5$ while the radius of the inclusions is 
		$R_I=1.25$. }\label{fig:2DStress_undef}
\end{figure}

\begin{figure}[h!]
\centering
\subfloat[]{\includegraphics[width=0.3\textwidth]{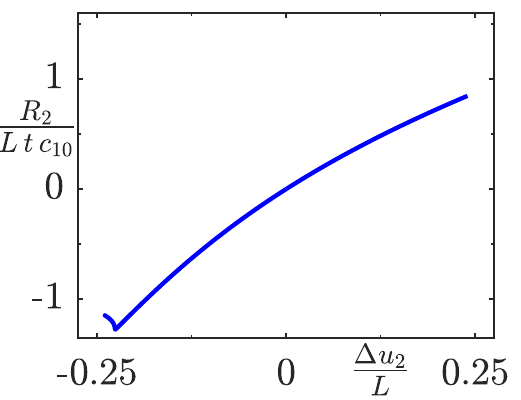}\label{fig:Pattern2D03_a}}\qquad
\subfloat[]{\includegraphics[width=0.3\textwidth]{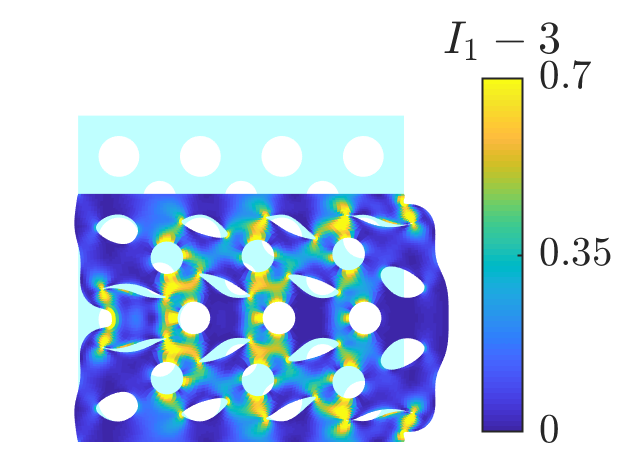}\label{fig:Pattern2D03_b}}\qquad
\subfloat[]{\includegraphics[width=0.3\textwidth]{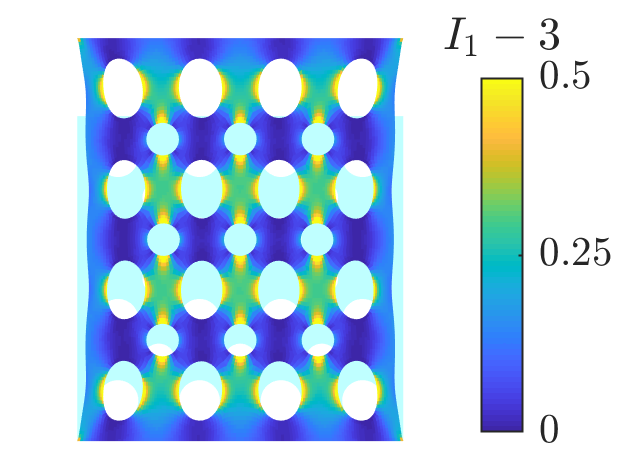}\label{fig:Pattern2D03_c}}
\caption{Plane stress problem, simulation results. (a) Normalized force-displacement plot, 
	$\Delta u_2$ is the vertical displacement applied to the top side of the boundary, $R_2$ is the 
	total reaction force conjugated to $\Delta u_2$, $L$ is the side of the domain, $t$ is the 
	thickness.
	Deformed configuration at maximum displacement in compression (b) and tension (c), $c_{10}$ is the modulus of the 	Neo-Hookean material model, the colormap is based on the first invariant of the deformation tensor.}\label{fig:Pattern2D}
\end{figure}

In this section we consider the equilibrium of an hyperelastic plate subjected to in-plane boundary conditions.
In particular we consider the domain, and the boundary conditions shown in figure \ref{fig:2DStress_undef}, where a vertical displacement $\Delta u_2$ is applied to the points of the top boundary, while the displacements of the points of the bottom boundary are fixed.
The domain features two types of internal boundaries.
The internal boundaries with radii $R_0$ are empty cavities, where standard traction free boundary conditions applied.
The internal boundaries with radii $R_I$ are rigid, circular, frictionless inclusions, whose points are constrained to remain at a fixed distance from the centre of the inclusion, which is free to move in both directions.
The material of the domain is assumed to be a Mooney-Rivlin type the modulus $c_{10}=10$ and $G=10^3$.
The presence of the inclusions has been introduced through Lagrange multipliers, by constraining the nodes lying on each inclusion boundary to remain on a circle having radius $R_I$, whose centre's coordinates were an unknowns, introduced through additional boundary conditions, determined at the equilibrium.
Figure \ref{fig:Pattern2D}.a shows the plot of the normalized constraint reaction as a function of the displacement, while figure \ref{fig:Pattern2D}.a shows the equilibrium configuration. 
As it can be observed, while the cavities dramatically changed their shape, both in compression and in tension, the inclusions maintained their circular shape.
The solution was achieved in 150 increments for the compressive branch and 100 increments for the tensile branch, with each increment converging in 4 or 5 iterations.

\subsection{Cylindrical symmetry problem with internal volume constraint}\label{sec:cylindric_symmetry}
In this section we discuss the solution of problems with cylindrical symmetry.
In the cases where the geometry of the domain, the material, and boundary conditions all have cylindrical symmetry, the problem can be significantly simplified incorporating the symmetry conditions within the displacement model.
Assuming that the axis of symmetry coincides with the $\bm{\imath}_2$ axis, in the absence of torsion, the deformation gradient takes the following form
\begin{equation}\label{F_cylindrical_symmetry}
\bm{F}^{\text{cyl}} = \left[\begin{matrix}
1+u_{1,1} & u_{1,2} & 0 \\
u_{2,1} & 1+u_{2,2} & 0 \\
0 & 0 & 1 + \frac{u_1}{X_1}
\end{matrix}\right] \,,
\end{equation}
and it is invariant to rotations around the symmetry axis.\\

\begin{figure}[h!]
	\centering
	\subfloat[]{\includegraphics[height=4cm]{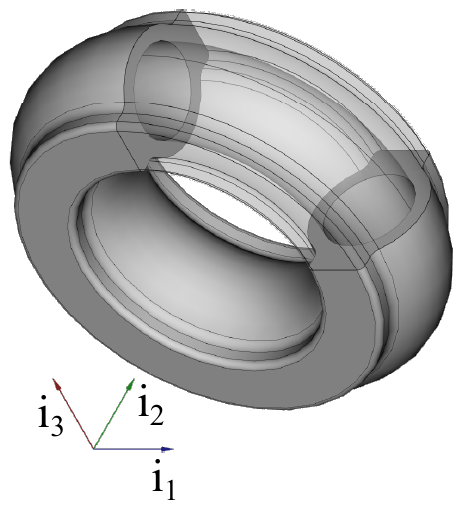}\label{fig:AxSymDomain3D04_02}}\hspace{2cm}
	\subfloat[]{\includegraphics[height=5cm]{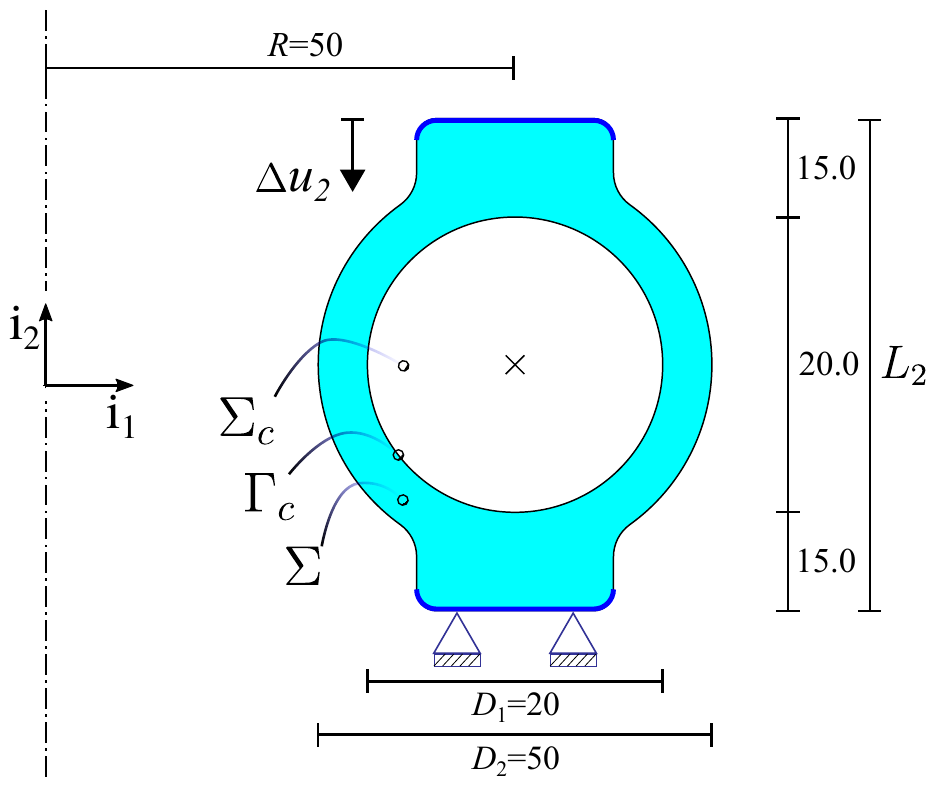}\label{fig:AxSymDomain04_02}}
	\caption{Axi-symmetric problem, domain geometry, dimensions and boundary conditions, $R$ is the distance of the centre of the section from the symmetry axes, the fillet radius are 2.0 and 5.0. The thicker line in (b) marks the portion of the boundary where displacement constraints were applied.}\label{fig:AxSymDomain04_02_all}
\end{figure}

In this example we consider the domain shown in figure \ref{fig:AxSymDomain3D04_02}, obtained by completely a complete rotation of the section shown in figure \ref{fig:AxSymDomain04_02} around the axis $\bm{\imath}_{2}$.
We assume that the material follows a Neo-Hookean model with modulus $c_{10}=10.0$, as per equation \eqref{NeoHookean}.
In addition to the boundary conditions illustrated in figure \ref{fig:AxSymDomain04_02}, we consider two cases for the behaviour of the internal cavity.
In one case we assume that the cavity is filled with an incompressible fluids.
For this case no particular shape is enforced, but only the value of the cavity's volume is kept constant during the solution.
In the second case we consider the cavity as an empty volume that can take any shape, with no other constraint than the external boundary conditions.\\

\begin{figure}[h!]
	\centering
	\subfloat[]{\includegraphics[height=4cm]{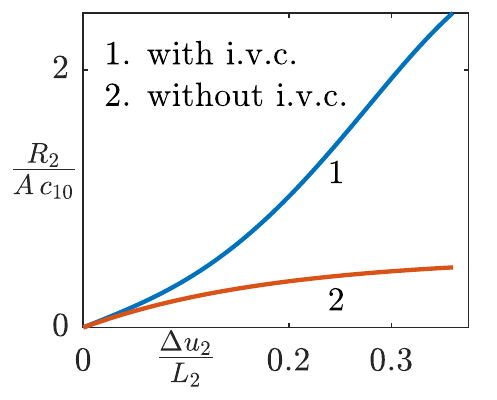}\label{fig:AxSymDomain04_02NH_a}}
	\subfloat[]{\includegraphics[height=4cm]{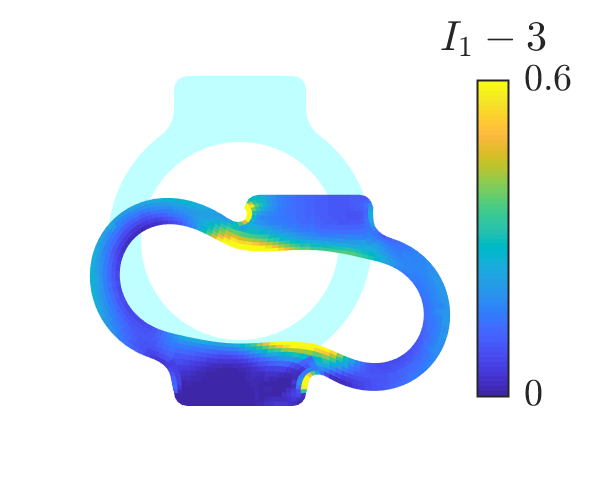}\label{fig:AxSymDomain04_02NH_b}}
	\subfloat[]{\includegraphics[height=4cm]{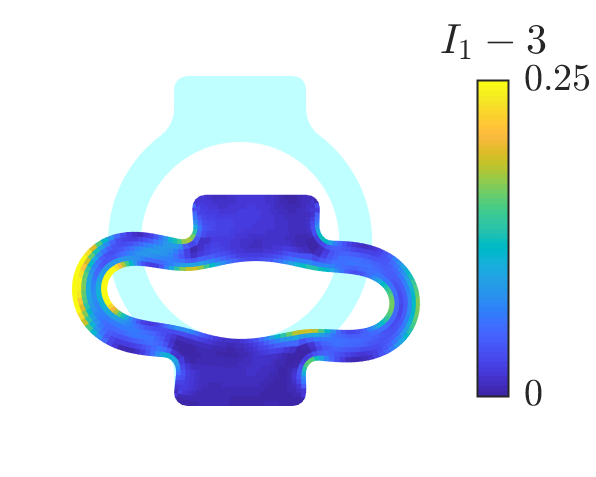}\label{fig:AxSymDomain04_02NH_c}}
	\caption{Axi-symmetric problem, simulation results. (a) Normalized force displacement plot, $\Delta u_2$ is the displacement applied in direction 2, $L_2$ is the initial height of the domain, $R_2$ is the total reaction force, $A$ is the area over which the boundary condition is applied, $c_{10}$ is the modulus of the Neo-Hookean material. Line 1. is the response with internal volume constraint (i.v.c.), line 2. is the response without i.v.c. . (b) Deformed cross section with i.v.c., and (c) without i.v.c. at maximum deformation.}\label{fig:AxSymDomain04_02NH}
\end{figure}

For any given configuration of the body, the volume of the internal cavity is given as
\begin{equation}\label{V_c0}
V_{c} = 2\pi \int_{\Sigma_c} x_1\, \mathrm{d}\Sigma = \pi 
\int_{\Gamma_c} x_1\,x_2\, \mathrm{d}x_1
\end{equation}
where $\Sigma_c$ is the intersection of the internal cavity volume with the plane $\bm{\imath}_{1}-\bm{\imath}_{2}$ in the initial configuration, and $\Gamma_c$ its boundary, as shown in figure \ref{fig:AxSymDomain04_02}. 
We recall that the second equality in equation \eqref{V_c0} holds thanks to the Gauss-Green theorem, and allows replacing the area integral with a curvilinear integral, without the need of discretizing the interior of the cavity.
The constraint on the volume of the inner cavity can then be introduced by means of Lagrange multipliers by requiring that the following holds
\begin{equation}\label{V_c0_consttraint}
V_{c}\left(u_k\right) = V_{c_0} \,,
\end{equation}
where $V_{c}\left(u_k\right)$ is the current volume of the cavity and  $V_{c_0}$ is the initial volume.
Therefore, the expression of the Lagrange functional to be minimized in order to solve the problem is the following
\begin{equation}\label{L_sym}
L^{\text{cyl}}\left(u_k, \lambda_c\right) = \Phi^{\text{cyl}}(u_k) - \lambda_c \left[V_{c}\left(u_k\right)-V_{c_0}\right]\,,
\end{equation}
where $\Phi^{\text{cyl}}$ is the deformation energy of the body, given as 
\begin{equation}\label{Phi_2D_stress_cyl}
\Phi^{\text{cyl}} = 2\pi \int_{\Sigma} \phi\left(\bm{F}^{\text{cyl}}\right) \, X_1 \mathrm{d}\Sigma\,,
\end{equation}
where $\Sigma$ is the section of the domain on the plane $\bm{\imath}_{1}-\bm{\imath}_{2}$, $X_1$, the first coordinate, is the distance from the rotation axis, and $\phi$ is the strain energy density of the material.
Figure \ref{fig:AxSymDomain04_02NH} shows the results of the simulations for the axisymmetric problem for a Neo-Hookean material model with modulus $c_{10}=10$. 
As expected, the presence of an incompressible fluid within the cavity, introduced through the constraint \eqref{V_c0_consttraint}, results in a general macroscopic stiffening of the solid, which enforces a lager widespread of the deformation across the domain, and lager average value of the deformation.
In both cases the solution was achieved in 200 steps with each steps taking 5 iterations to converge in the case with the cavity volume constraint, and 6 iterations in the cases without the constraint.
\\

We remark that the cylindrical symmetry, in the above example, was introduced by simply incorporating it in the displacement model and in the definition of $\bm{F}$, in equation \eqref{F_cylindrical_symmetry}, while at no point, in the statement of the elastic problem, it was necessary to express the equilibrium equation in cylindrical coordinates.

\subsection{Three-dimensional solid with large geometrical non-linearities}\label{sec:3D_solid}

\begin{figure}[h!]
\centering
\subfloat[]{\includegraphics[height=4cm]{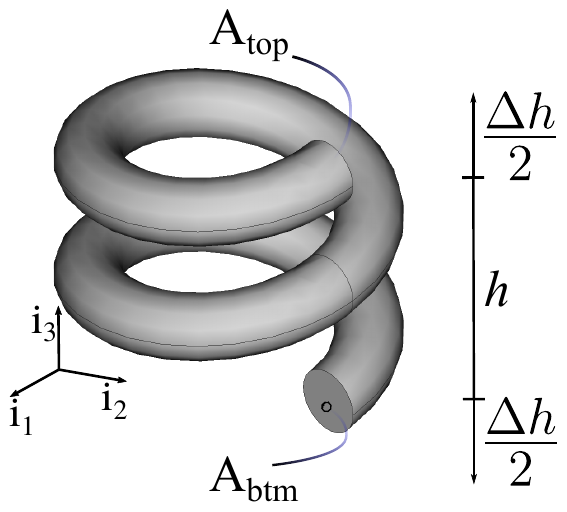}\label{fig:3DSpringCAD}}\qquad
\subfloat[]{\includegraphics[height=4cm]{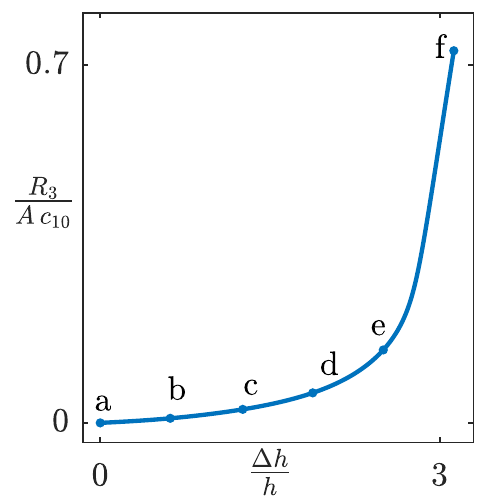}\label{fig:3DSpringFineMeshNHb}}\qquad
\subfloat[]{\includegraphics[height=4cm]{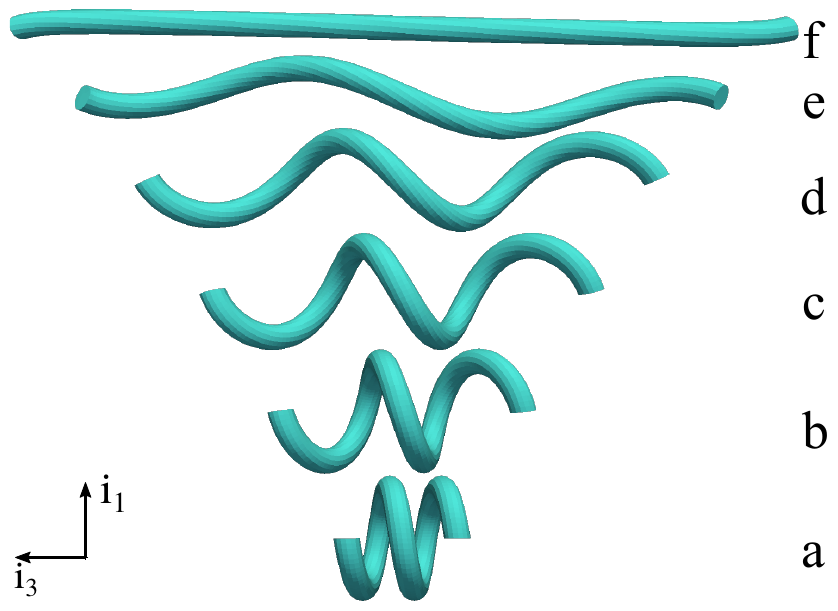}\label{fig:SpringFineMeshNHb}}
\caption{Three-dimensional solid. (a) Domain's geometry, $A_{\text{top}}$ and $A_{\text{btm}}$ are the top and the bottom end cross sections of the helix, respectively; (b) Reaction force,  $R_3$, $A$ is the cross section area, versus normalized applied displacement, $c_{10}$ is the modulus of the Neo-Hookean model; (c) Deformed configurations at the stages of the simulation marked with dots in (b).} \label{fig:3DSpring}
\end{figure}

In this section we analyse the response of a three dimensional hyperelastic solid undergoing large displacements.
Figure \ref{fig:3DSpringCAD} shows the domain geometry, which consists in a right-handed helicoidal solid with circular cross section.
The helix has radius $R_e=20$, pitch $p=20$ and height $h=40$.
In this example the boundary conditions have been applied constraining the displacement of the centre of mass of the two ends of the helix to move along direction $\bm{\imath}_{3}$, increasing the height of the helix, has shown in figure \ref{fig:3DSpringCAD}.
Therefore the following set of equation was imposed on the nodes of the end cross sections 
\begin{subequations}
\begin{align}
\int_{A_{\text{btm}}} u_1\, \mathrm{d}A &=0\qquad,\qquad \int_{A_{\text{top}}} u_1\, \mathrm{d}A= 0 \label{eq_spring_a}\\
\int_{A_{\text{btm}}} u_2\, \mathrm{d}A &=0\qquad,\qquad \int_{A_{\text{top}}} u_2\, \mathrm{d}A= 0 \label{eq_spring_b}\\
\int_{A_{\text{btm}}} u_3\, \mathrm{d}A &=-\frac{\Delta h}{2}\quad,\quad \int_{A_{\text{top}}} u_3\, \mathrm{d}A=\frac{\Delta h}{2} \label{eq_spring_c}
\end{align}\label{eq_spring}
\end{subequations}
where $A_{\text{top}}$ and $A_{\text{btm}}$ are the top and the bottom cross section, respectively.
The boundary conditions have been applied, similarly to the previous example, through Lagrange multipliers, adding a constraint equation for each of the equations \ref{eq_spring}.
We observe that the value of Lagrange multipliers conjugated to equations \eqref{eq_spring_c}, at equilibrium correspond to the total constraint reactions in direction $\bm{\imath}_3$ on the bottom and the top faces, respectively, which are the active forces, producing the deformation of the spring.
Because of the symmetry of the domain and of the boundary conditions, at equilibrium both the active reaction forces have the same value, $R_3$.
Figure \ref{fig:3DSpringFineMeshNHb} shows the normalized plot of $R_3$ as a function of $\Delta h$, and figure \ref{fig:SpringFineMeshNHb} shows the deformed configurations of the helix corresponding to the points marked in figure \ref{fig:3DSpringFineMeshNHb}.
As we can observe the plot of the reaction force shows the expected hardening behaviour due to the alignment of the helix with the applied force.
The solution was obtained in 300 increments, with each increment taking 3 iterations to converge.

\section{Concluding remarks}

Automatic differentiation (AD) techniques allow for the accurate and efficient numerical evaluation of the derivatives of a multivariate function.
In this paper, AD has been used for stating and solving non-linear finite element solid mechanics problems.
The approach presented here focuses in particular on Green elastic materials, for which the deformation work is an exact differential, and the solid can be treated as a proper conservative thermodynamic system.
In these cases, the residual force vector and the tangent stiffness matrix of the model coincide, respectively, with the gradient and the Hessian of the system's free energy, which can be numerically evaluated through AD.
Therefore, with the approach presented here no explicit calculation of the stress tensor, nor of the elasticity tensor is required, nor it is necessary to implement the complex kinematics that link the degrees of freedom of the model to the internal forces and their derivatives.
The same framework can also be applied with arbitrary, non conservative, material models, although, here the explicit calculation of the components of the stress tensor is required, while the calculation of the elasticity tensor and of the tangent stiffness matrix can still be automated.
In the same way, sophisticated constraints equations and boundary conditions, can be introduced by means of Lagrange multipliers, and treated through AD.
The method has been presented with a number of examples that illustrate the solution of selected non-linear problems, featuring hyperelastic material models, and complex constraints, along with the computer programs developed for producing the results included in this article.



\bibliography{bibliography}

\pagebreak 
\appendix
\section*{Appendices}
\section{Implementation of dual number systems in the Julia programming language}\label{sec:implementation_in_Julia}

In this section we illustrate a possible implementation of the dual number system in the Julia programming language.
Julia is a dynamically typed scientific programming language, whose semantic is particularly suitable for the description of physical problems \citep{Bezanson2017,Perkel2019}.
Aside to user-defined data types, Julia permits to overload existing operators or functions to be evaluated on the new types.
Therefore the same script that evaluates a numeric function on floating point numbers, can be used to operate on dual numbers, once their arithmetic has been implemented, and produce dual number as a result.\\

The script blocks reported below show a possible implementation of dual numbers in Julia.
The dual number type is called D2 and it is defined in Listing \ref{src:dual_number_definition} as having a scalar component  {\footnotesize\ttfamily v}, that stores the current value of the variable, a one dimensional array, {\footnotesize\ttfamily d1}, that stores all of the first derivatives of {\footnotesize\ttfamily v}, and a two dimensional array, {\footnotesize\ttfamily d2}, that stores all of the second derivatives of {\footnotesize\ttfamily v}.
\begin{lstlisting}[caption={Definition of dual numbers in Julia}, label={src:dual_number_definition}]
struct D2 <: Number
	v::Float64
	d1::Array{Float64,1}
	d2::Array{Float64,2}
end
\end{lstlisting} 
The code in Listing \ref{src:dual_operators_definition} extends some ordinary maths operators to function with the \texttt{D2} type.
The first line in the script block informs the language that the scope of the mentioned operators, defined in the Base module, will be extended, while the remaining lines implement the arithmetic of dual numbers as defined in section \ref{sec:dual_number_field}, where each component of a dual number value is accessed through the dot syntax ({\footnotesize\ttfamily .}), and the single quote ({\footnotesize\ttfamily\textquotesingle}) denotes array transposition.
\begin{lstlisting}[caption={Operators overloading}, label={src:dual_operators_definition}]
import Base: +,-,*,/,^
+(x::D2, y::D2) = D2(x.v+y.v, x.d1+y.d1, 
		x.d2+y.d2)
-(x::D2, y::D2) = D2(x.v-y.v, x.d1-y.d1,
		x.d2-y.d2)
*(x::D2, y::D2) = D2(x.v*y.v, 
	x.d1*y.v+y.d1*x.v,
	x.d2*y.v+x.v*y.d2+
	x.d1*y.d1'+y.d1*x.d1')
/(x::D2, y::D2) = D2(x.v/y.v,
	x.d1/y.v-(x.v/y.v^2)*y.d1,
	x.d2/y.v-(x.d1*y.d1'+y.d1*x.d1')/
	y.v^2+2x.v*(y.d1*y.d1')/y.v^3-
	(x.v/y.v^2)*y.d2)
^(x::D2, n::Int64) = D2(x.v^n,
	(n*x.v^(n-1))*x.d1,
	(n*(n-1)*x.v^(n-2))*
	(x.d1*x.d1')+(n*x.v^(n-1))*x.d2)
\end{lstlisting}
In Listing \ref{src:dual_number_example_float} the function given by equation \eqref{y(x1,x2,x3)} is defined in the first line, in the following line it is evaluated for the double precision floating point values \texttt{x1=x2=x3=1.0} and the result is printed.
\begin{lstlisting}[caption={Numerical example with \texttt{Float64} arguments}, label={src:dual_number_example_float}]
y(x1,x2,x3) = x1^3*x2^2 + x3^2
println("\n y0: ", y(1., 1., 1.))

 y0: 2.0
\end{lstlisting} 
In Listing \ref{src:dual_number_example}, the variables \texttt{x1}, \texttt{x2} and \texttt{x3} are defined as dual quantities, of the \texttt{D2} type, where the first argument of the call to the \texttt{D2} constructor is the value of the variable, the second argument is the gradient of each variable, and the third argument is the Hessian.
We remark that independent variables are defined by properly specifying the components of their gradient and Hessian.
In fact, independent variables are such that the derivative with respect to themselves is one, while all other derivatives are nought.
Thus \texttt{x1} is the independent variable that occupies the first position of the gradient, \texttt{x2} the second and \texttt{x3} the third.
As a result, the derivative with respect to \texttt{x1} of any operation involving \texttt{x1}, \texttt{x2} and \texttt{x3}, with will be stored in the first component of the \texttt{d1} array of the result, and similarly for the derivatives with respect to \texttt{x2} and \texttt{x3}, and for higher order derivatives.
Following in Listing \ref{src:dual_number_example}, the same function y, defined in Listing \ref{src:dual_number_example_float}, is called with the dual quantities just defined and the result, which is a dual quantity itself, is printed.
As we can observe, the result returned by the function this time includes both the value of the  function \texttt{y(x1,x2,x3)}, and its gradient and Hessian.
\begin{lstlisting}[caption={Numerical example with \texttt{D2} arguments}, label={src:dual_number_example}]
x1 = D2(1., [1., 0, 0], zeros(3,3))
x2 = D2(1., [0, 1., 0], zeros(3,3))
x3 = D2(1., [0, 0, 1.], zeros(3,3))
println("\n yd: ", y(x1, x2, x3))

 yd: D2(2.0,[3.0, 2.0, 2.0],
 	[6.0 6.0 0.0; 
	 6.0 2.0 0.0; 
	 0.0 0.0 2.0])
\end{lstlisting} 
We observe that having overloaded the operators involved in the definition of \texttt{y(x1,x2,x3)} allowed us to call the same function with both data type without making any modification or having to add any specification to the function itself.\\

We remark that the implementation of dual numbers in Julia as presented in this section is an attempt to provide a brief and clear illustration of a possible computer implementation of AD, through operators overloading, nonetheless in this form it does not exploit any of the  powerful features offered by the Julia programming language, like parametric types, and macros \citep{Perkel2019,Bezanson2017}.
The implementation developed for the solution of the example presented in the paper, available through \citep{Mendeley}, which is based on \cite{RevelsLubinPapamarkou2016}, makes a better use of Julia's features and functionalities and ensures better performances than the example presented in this section.

\section{Arbitrary order dual number systems}\label{sec:arbitrary_order_dual_numbers}
In this section we briefly generalize the definition of dual numbers to an arbitrary order of differentiation.
Let $\bm{x}$ be a dual number of dimension $N$ and order $K$
\begin{equation}\label{mathbfx_ii}
\begin{split}
\bm{x} \equiv x_0 &+ x_{i_1} \imath_{i_1} + x_{i_1i_2} \,\imath_{i_1i_2} + x_{i_1i_2i_3} 
\,\imath_{i_1i_2i_3} + \cdots \\ 
& +  x_{i_1\dots i_K} \,\imath_{i_1\dots i_K} \quad \text{with} \qquad \begin{cases}
i_1 &\in 1\dots N \\
i_2 &\in i_1 \dots N \\
i_3 &\in i_2 \dots N \\
&\vdots \\
i_{K} &\in i_{K-1} \dots N \\
\end{cases}
\end{split}\end{equation}
with $\imath_{j}$ the canonical base of $\mathcal{R}^N$, with $j \in 1 \dots N$, and $\bm{\imath}_i, \bm{\imath}_{ij}, \imath_{ijk} \dots, \imath_{i_1\dots i_K}$ are symbols defined as
\begin{equation}\label{imath_ij_ii}
\begin{aligned}
\imath_{i_1 i_2} &\equiv \imath_{i_1}\otimes \imath_{i_2} + \imath_{i_2}\otimes \imath_{i_1}\\
\imath_{i_1 i_2 i_3} & \equiv \imath_{i_1} \otimes \imath_{i_2} \otimes \imath_{i_3} + \imath_{i_1} \otimes \imath_{i_3}\otimes \imath_{i_2} + \imath_{i_3} \otimes \imath_{i_1}\otimes \imath_{i_2} + \\   & \hspace{20pt} \imath_{i_3} \otimes \imath_{i_2}\otimes \imath_{i_1} + \imath_{i_2} \otimes \imath_{i_3}\otimes \imath_{i_1}  + \imath_{i_2} \otimes \imath_{i_1}\otimes \imath_{i_3} \\ 
&  \hspace{20pt} \vdots \hspace{100pt} \vdots \\
\imath_{i_1\dots i_K} & \equiv \sum_{I^K \in \Pi(K)} \imath_{I^K_1} \otimes \cdots \otimes \imath_{I^K_K}
\end{aligned}
\end{equation}
where $I^K$ is a permutation of the indices $1 \dots K$, $I_i^K$ are its elements, and $\Pi(K)$ is the set of all the permutations of the indices $1 \dots K$. With respect to equations \eqref{imath_ij_ii} we observe that the following holds
\begin{equation}\label{i_ij=i_ji_ii}
\begin{aligned}
\bm{\imath}_{ij} &= \bm{\imath}_{ji}\\
\imath_{ijk} &= \imath_{ikj} = \imath_{jik} = \imath_{jki} = \imath_{kij} =\imath_{kji} \\
&  \hspace{20pt} \vdots \hspace{100pt} \vdots \\
\imath_{I^K} &= \imath_{J^K} \hspace{50pt} \forall \, I^K,J^K \in  \Pi(K)
\end{aligned}.
\end{equation}
The quantities $x_0, x_i, x_{ij}, x_{ijk} \dots, x_{i_1\dots i_K}$ are real scalars and are the 
components of $\bm{x}$, $x_0$ is the real part of $\bm{x}$, the remaining are dual parts of order 
$1, 2, \dots K$.
Two dual numbers of dimension $N$ and order $K$ are identical if and only if all of their components are identical, as follows
\begin{equation}\label{x=y}
\bm{x}=\bm{y} \quad \iff \quad \begin{cases}
x_0 = y_0 \\
x_i = y_i \\
\vdots \qquad \vdots \\
x_{i_1\dots i_K} = y_{i_1\dots i_K}
\end{cases}
\end{equation}
The sum of two dual numbers is defined as the dual number whose components are the sum of the components, as follows
\begin{equation}\label{z=x+y}
\bm{z} = \bm{x}+\bm{y} \quad \iff \quad \begin{cases}
z_0 = x_0 + y_0 \\
z_i = x_i + y_i \\
\vdots \hspace{50pt} \vdots \\
z_{i_1\dots i_K} = x_{i_1\dots i_K} + y_{i_1\dots i_K}
\end{cases}
\end{equation}
The product of two dual  numbers is a dual number obtained as the sum of the mixed products of their components, where the following rules applies for the product of the symbols $\bm{\imath}_i, \bm{\imath}_{ij}, \imath_{ijk} \dots, \imath_{i_1\dots i_K}$
\begin{equation}\label{i_i i_j i_k = i_ijk}
\begin{aligned}
\bm{\imath}_i \bm{\imath}_j &\equiv \bm{\imath}_{ij} \\
\bm{\imath}_i \bm{\imath}_j \imath_k &= \imath_{i} \imath_{jk} \equiv \imath_{ijk}\\
& \vdots \hspace{50pt} \vdots \\
\bm{\imath}_{1} \dots \imath_K  &= \bm{\imath}_{1} \imath_{2\dots K} \equiv \imath_{1\dots K}
\end{aligned}\qquad ,
\end{equation}
\begin{equation}\label{i_li_ijk=0}
\bm{\imath}_i\imath_{1\dots K} \equiv 0 \quad ,
\end{equation}
where equations \eqref{i_i i_j i_k = i_ijk} produce the contribution to higher terms in the product as results of the products of lower order terms in the factors, and equation \eqref{i_li_ijk=0} ensures that no component with order higher than $K$ appears in the result.
The components of the product are given as
\begin{equation}\label{z=xy_ii}
\bm{z} = \bm{x} \bm{y} \iff \begin{cases}
z_0 = x_0  y_0 \\
z_i = x_iy_0 + x_0y_i\\
z_{ij} = x_{ij}y_0 + x_i y_j + x_j y_i + x_0 y_{ij} \\
z_{ijk} = x_{ijk}y_0 + x_{ij} y_k + x_i y_{jk} + x_0 y_{ijk} \\
\quad \vdots \hspace{100pt} \vdots \\
z_{i_1\dots i_K} = x_{i_1\dots i_K}y_0 + x_{i_1\dots i_{K-1}}y_{i_K} + \dots + x_0y_{i_1\dots i_K}
\end{cases},
\end{equation}
With reference to the quotient of two dual numbers, we observe that this operation is equivalent to the product of the first time the inverse of the second, where the inverse of a dual number is obtained by solving the following
\begin{equation}\label{1/x}
\frac{1}{\bm{x}} = \bm{y} \iff \ \bm{y}\bm{x} = 1,
\end{equation}
from which it results
\begin{equation}\label{1/x_ii}
\frac{1}{\bm{x}} = \bm{y} \iff \begin{cases}
y_0x_0 = 1 \\
x_iy_0 + x_0y_i =0 \\
x_{ij}y_0 + x_i y_j + x_j y_i + x_0 y_{ij} =0\\
x_{ijk}y_0 + x_{ij} y_k + x_i y_{jk} + x_0 y_{ijk} =0\\
\quad \vdots \hspace{100pt} \vdots \\
x_{i_1\dots i_K}y_0 + x_{i_1\dots i_{K-1}}y_{i_K} + \dots + x_0y_{i_1\dots i_K} =0 
\end{cases},
\end{equation}
where we observe that the right hand side of the definition \eqref{1/x_ii} is an lower diagonal linear system in the unknowns $y_{\cdots}$, which can be easily solved by back-substitution starting from the first equation.

\end{document}